%% file: ejs-BraultKeribinMariadassou.tex
\documentclass[ejs]{imsart}

\RequirePackage[OT1]{fontenc}
\RequirePackage{amsthm,amsmath,amssymb,amsfonts,dsfont,mathtools}

\RequirePackage[colorlinks,citecolor=blue,urlcolor=blue]{hyperref}

\usepackage{tikz}
\usepackage{longtable}
\usepackage{version}
\usepackage{textcomp}
\usepackage{url}

\usepackage[babel=true]{csquotes}

\newcommand{\CK}[1]{\textcolor{red}{ #1}}
\newcommand{\VB}[1]{\textcolor{purple}{#1}}
\newcommand{\MM}[1]{\textcolor{orange}{#1}}
\input{commandes.tex}

\startlocaldefs
\numberwithin{equation}{section}
\theoremstyle{plain}
\newtheorem{thm}{Theorem}[section]
\newtheorem{dof}[thm]{Definition}
\newtheorem{proposition}[thm]{Proposition}
\newtheorem{corollaire}[thm]{Corollary}
\newtheorem{lemme}[thm]{Lemma}
\theoremstyle{remark}
\newtheorem{rem}[thm]{Remark}

\endlocaldefs

\doi{10.1214/154957804100000000}
\pubyear{0000}
\volume{0}
\firstpage{0}
\lastpage{0}
\arxiv{1704.06629}

\startlocaldefs
\endlocaldefs

\begin{document}

\begin{frontmatter}

% "Title of the Paper"
\title{Consistency and Asymptotic Normality of Latent Block Model Estimators}
%\thankstext{t1}{This is an original survey paper}
\runtitle{Asymptotics for LBM}

% indicate corresponding author with \corref{}
% \author{\fnms{John} \snm{Smith}\thanksref{t2}\corref{}\ead[label=e1]{smith@foo.com}\ead[label=e2,url]{www.foo.com}}
% \thankstext{t2}{Thanks to somebody} 
% \address{line 1\\ line 2\\ \printead{e1}\\ \printead{e2}}

\author{\fnms{Vincent} \snm{Brault}\ead[label=e1]{vincent.brault@univ-grenoble-alpes.fr}\ead[label=e2,url]{https://www-ljk.imag.fr/membres/Vincent.Brault/}}
\address{Univ. Grenoble Alpes, CNRS, Grenoble INP\footnote{Institute of Engineering Univ. Grenoble Alpes}, LJK, 38000 Grenoble, France\\ \printead{e1}\\ \printead{e2}}
\and
\author{\fnms{Christine} \snm{Keribin}\ead[label=e3]{christine.keribin@math.u-psud.fr}\ead[label=e4,url]{https://www.math.u-psud.fr/\~{ }keribin/}}
\address{Universit\'e Paris-Saclay, CNRS, Inria, Laboratoire de math\'ematiques d'Orsay, 91405, Orsay, France\\ \printead{e3}\\ \printead{e4}}
\and
\author{\fnms{Mahendra} \snm{Mariadassou}\ead[label=e5]{mahendra.mariadassou@inrae.fr}\ead[label=e6,url]{https://mahendra-mariadassou.github.io/}}
\address{MaIAGE, INRAE, Universit\'e Paris-Saclay, 78352 Jouy-en-Josas, France\\ \printead{e5}\\ \printead{e6}}

\runauthor{V. Brault, C. Keribin and M. Mariadassou}

\begin{abstract}
The Latent Block Model (LBM) is a model-based method to cluster simultaneously the 
$\dd$ columns and $\n$ rows of a data matrix. Parameter estimation in LBM is a 
difficult and multifaceted problem. Although various estimation strategies have 
been proposed and are now well understood empirically, theoretical guarantees 
about their asymptotic behavior is rather sparse and most results are 
limited to the binary setting. We prove here theoretical guarantees in the 
valued settings. We show that under some mild conditions on the 
parameter space, and in an asymptotic regime where $\log(\dd)/\n$ and 
$\log(\n)/\dd$ tend to $0$ when $\n$ and $\dd$ tend to infinity, (1) the 
maximum-likelihood estimate of the complete model (with known labels) is 
consistent and (2) the log-likelihood ratios are equivalent under the complete 
and observed (with unknown labels) models. This equivalence allows us to 
transfer the asymptotic consistency, and under mild conditions, asymptotic 
normality, to the maximum likelihood estimate under the observed model. 
Moreover, the variational estimator is also consistent and, under the same 
conditions, asymptotically normal. 
%These results can be used to show that the integrated complete likelihood (ICL) criterion has the same behavior as BIC under the true distribution.
\end{abstract}

%\begin{keyword}[class=MSC]
%\kwd[Primary ]{}
%\kwd{}
%\kwd[; secondary ]{}
%\end{keyword}


% history:
% \received{\smonth{1} \syear{0000}}

%\tableofcontents

\end{frontmatter}

%%%%%%%%%%%%%%%%%%%%%%%%%%%%%%%%%%%%%%%%%%%%%%%%%%%%%%%%%%%%%%%%%%%%%%
\section{Introduction}

\begin{figure}[ht!]
\centering
\includegraphics[scale=0.08]{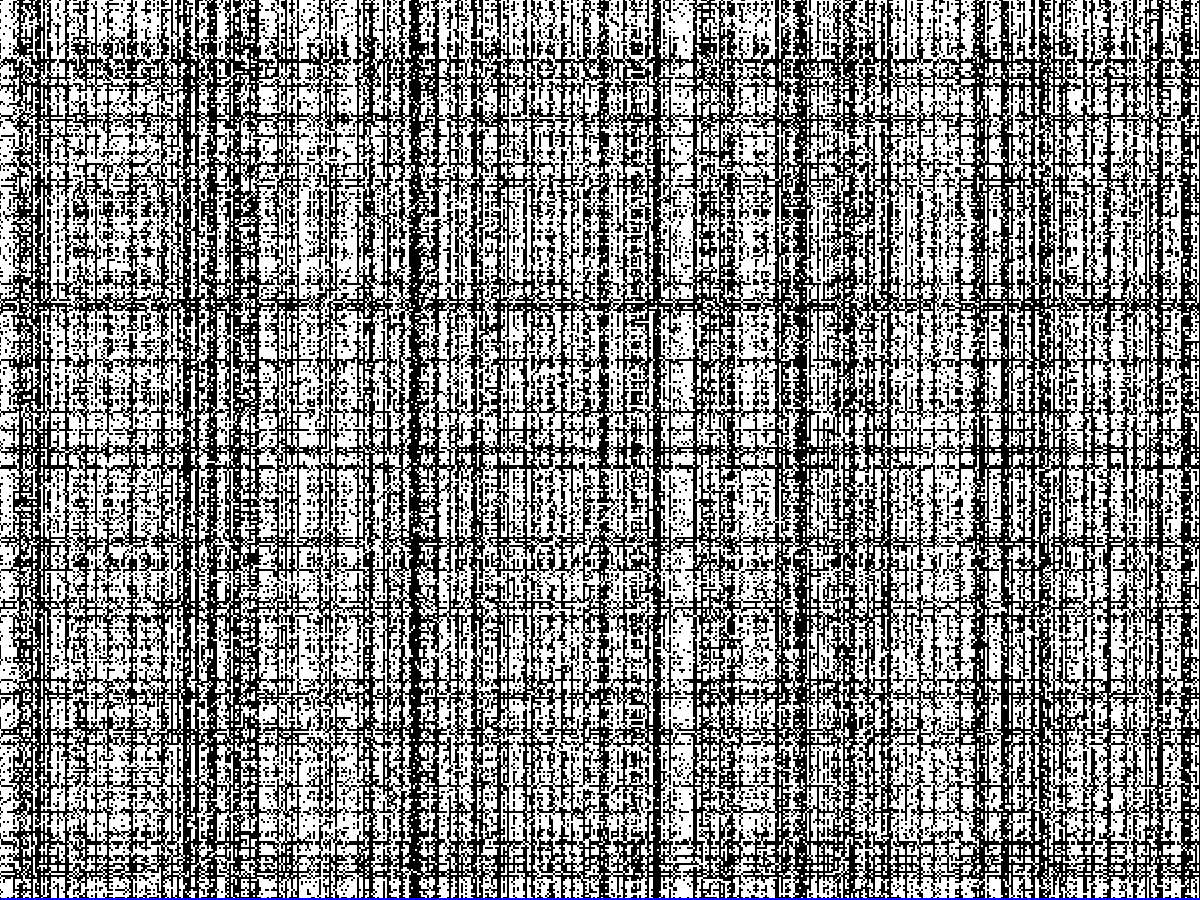}
\includegraphics[scale=0.08]{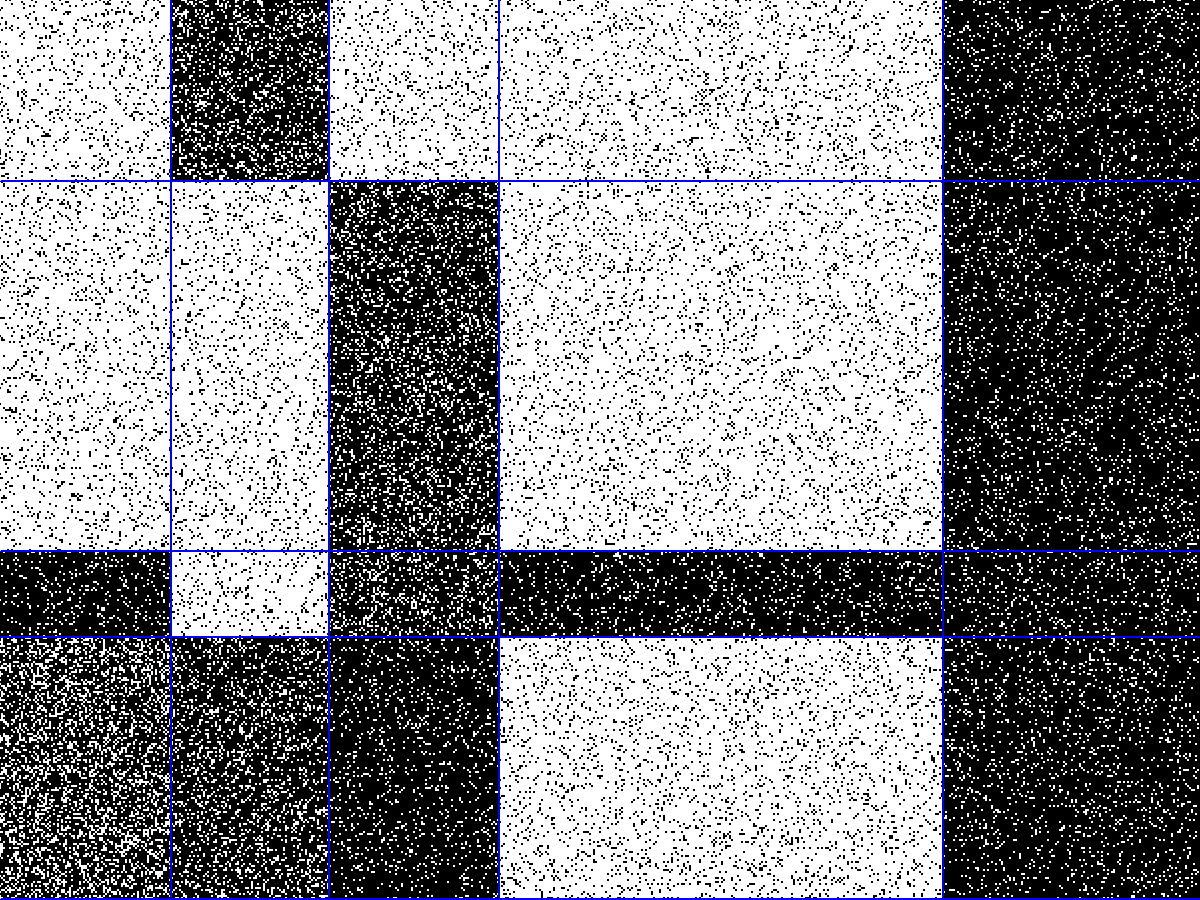}
\includegraphics[scale=0.08]{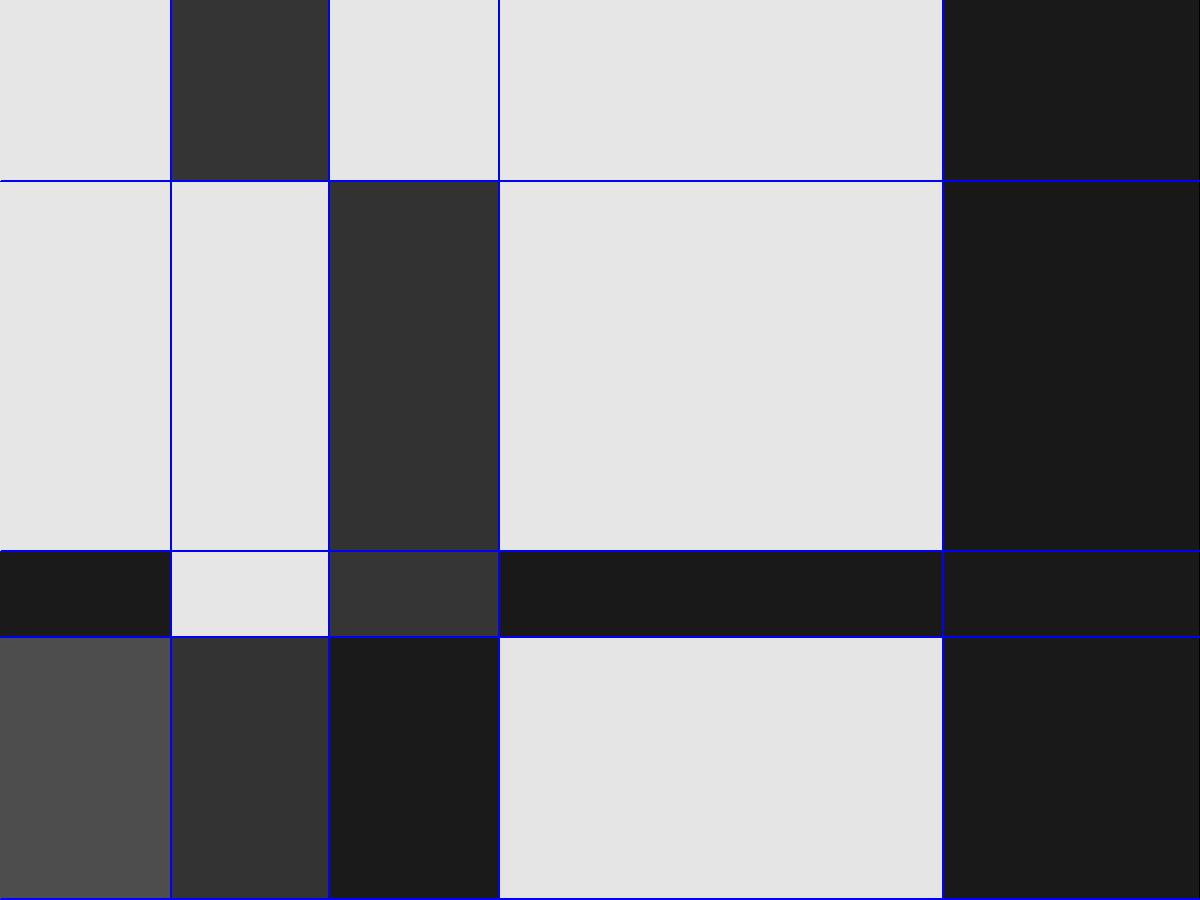}
\caption{A binary data matrix before (left) and after (middle) row and column reordering, and its parsimonious data representation (right).}\label{fig:cc}
\end{figure}

Co-clustering is an unsupervised method to cluster simultaneously the $\n$ rows and $\dd$ columns of a rectangular data matrix.  The assignments of each row to one of the row-clusters and of each column to one of the column-clusters are unknown and the aim is to determine them. Then, rows and columns can be re-ordered according to 
their assignments, highlighting the natural structure of the data with distinct blocks having homogeneous observations. This leads to a parsimonious data representation, as can be shown on Figure \ref{fig:cc}.

%\begin{figure}[ht!]
%\centering
%\includegraphics[scale=0.10]{figures/Tab1bis}
%\includegraphics[scale=0.10]{figures/Tab2bis}
%\includegraphics[scale=0.05]{figures/Tab3bis}
%\caption{A binary data matrix before (left) and after (middle) row and column reordering, and its parsimonious data representation (right).}\label{fig:cc}
%\end{figure}

Co-clustering can be used in numerous applications, and especially ones with large data sets, such as   recommendation systems (to discover a segmentation of customers with regard to a segmentation of products), genomics (to simultaneously define  
groups of genes having the same expression with regards to groups of 
experimental conditions) or text mining (to define simultaneously groups of 
texts and groups of words).

  Among the co-clustering methods, the Latent Block Model (LBM) defines a probabilistic model as a mixture model with latent rows and columns assignments.
  %is based on the definition of a probabilistic model.   
LBM can deal with binary (\cite{govaert2008block}), Gaussian 
(\cite{govaert2013co}), categorical (\cite{keribin2015estimation}) or count 
(\cite{govaert2010latent}) data. Due to the complex dependency structure induced by this modeling, 
neither the likelihood, nor the distribution of the 
assignments conditionally to the observations needed in the E-step of the EM algorithm, traditionnally used for mixture models, are  numerically tractable. Estimation can be however performed either with a variational approximation leading to an approximate value of the maximum likelihood estimator, or with a Bayesian approach 
(VBayes algorithm or Gibbs sampler). For example, \cite{keribin2015estimation} 
recommends using a Gibbs sampler combined with a VBayes algorithm. 

%These methods give satisfactory results in the applications, 
The asymptotics  of the maximum likelihood (MLE) and variational (VE) estimators also  raise interesting theoretical questions. % are still an open question.
This topic was first addressed for the stochastic block model (SBM)  (\cite{Snijders1997}), where the data is a random graph  encoded by its adjacency binary matrix: the rows and columns represent the nodes, so that there is only 
one partition, shared by rows and columns, and  a unique asymptotic direction.% as the number of rows and columns grow at the same rate.

%Some partial results exist for LBM, and this question has been solved for SBM (Stochastic Block Model), a special case of LBM where the data is a random graphe  encoded by its adjacency matrix (rows and columns represents the same units, so that there is only one partition, the same for rows and columns). 
For a binary SBM and under the true parameter value, Theorem 3 of \cite{celisse2012consistency} states that the distribution of the assignments conditionally to the 
observations converges to a Dirac of the real assignments. 
Moreover, this convergence remains valid under the estimated parameter value, 
assuming that this estimator converges at rate at least $\n^{-1}$, where $\n$ is 
the number of nodes (Proposition 3.8). This assumption is not trivial, and it 
was not established that such an estimator exists except in some particular cases 
(\cite{ambroise2012new} for example). \cite{mariadassou2015} presented a unified 
frame for LBM and SBM in case of valued observations satisfying a concentration 
inequality, and showed the consistency of the conditional distribution of the assignments under all parameter values in a neighborhood of the true value. 
\cite{bickel2009nonparametric} and \cite{bickel2013asymptotic} proved the 
consistency and asymptotic normality of the MLE for the binary SBM but 
failed to account for complications induced by symmetries in the parameter. 
Building upon the work from \cite{celisse2012consistency}, they first studied 
the asymptotic behavior of the MLE in the complete model (observations and 
assignments) with binary observations which is simple to handle; 
then, they showed that the complete likelihood and the marginal likelihood have 
similar asymptotic behaviors by the use of a Bernstein inequality for bounded 
observations.

%We extend these results to the double \CK{(rows, columns)}asymptotic  framework of LBM, following the way of \cite{bickel2013asymptotic}, and for observations coming from some exponential family. 
Following the main ideas of \cite{bickel2013asymptotic}, we prove that the observed likelihood ratio and the complete likelihood ratio computed at the 
true assignments are asymptotically equivalent, up to a multiplicative term. This term depends on some model symmetry and was omitted in \cite{bickel2013asymptotic} although it is necessary to prove the asymptotic results. We then settle the asymptotic normality of the maximum likelihood and variational estimators.
All these results are stated not only for binary observations, but also more 
generally for observations coming from univariate exponential families in canonical form, which is essential regarding the LBM usages. This leads us to develop a Bernstein-type inequality for sub-exponential variables as the Hoeffing's concentration 
inequality used in \cite{bickel2013asymptotic} is only relevant for 
upper-bounded observations. %The objective of this theoretical work is to prove the asymptotic behavior of Maximum Likelihood and Variational Estimators for LBM under the true model. We do not address here the choice of  the number of blocks. We refer to 
%In this work, we assume that the data are generated by a true model with known numbers of row-clusters and column-clusters.
%(which assumed the consistency) for the use of ICL in this case and to for other recent results in this topic.

%Moreover, we introduce the concept of model symmetry  which was not pointed out by these authors, but is necessary to set the asymptotic behavior. 
%The asymptotic normality of the variational estimator is also settled, and an application to model selection criteria is presented.

The paper is organized as follows. The model, main assumptions and notations are 
introduced in Section \ref {sec:ModelAndAssumptions}, where the concept of model 
symmetry is also discussed. Section \ref{sec:mle-complete-likelihood} 
proves the asymptotic normality of the complete likelihood estimator, and 
section \ref{sec:profile-likelihood} studies conditional and profile log-likelihoods. Our main result showing that the observed likelihood ratio
behaves like the complete likelihood ratio is stated in section 
\ref{sec:big-theorem}, and its consequences in terms of consistency and asymptotic normality of the MLE and variational estimators are presented in section \ref{sec:est-asympt-beh}. 
%We ended with an application to model selection and a discussion on sparsity.
Most of the proofs are postponed to the appendices to improve the general readibility :  appendix \ref{sec:proofPL} for properties of conditional and profile log-likelihoods,  \ref{sec:proofCVG} for the steps of the main result,  \ref{sec:subexp} for concentration inequalities for specific sub-exponential variables and  \ref{sec:technical-lemma} for other technical results. 
%settles three different types of assignment behaviors

%%%%%%%%%%%%%%%%%%%%%%%%%%%%%%%%%%%%%%%%%%%%%%%%%%%%%%%%%%%%%%%%%%
\section{Model, assumptions and definitions}
\label{sec:ModelAndAssumptions}
We observe a data matrix $\X=(\xij)$ with $\n$ rows and $\dd$ columns. The LBM assumes that there exists a latent structure in the form of the Cartesian product of a partition of $g$ row-clusters by a partition of $m$ column-clusters with the following characteristics:
\begin{itemize}
\item  the latent row assignments $\bz=(\bz_1,\ldots,\bz_n)$ are independent and identically distributed with a common multinomial distribution on $g$ categories:  \[\bz_i\sim \mathcal{M}(1, \bpi=(\pi_1,\ldots,\pi_g))\]
For $k=1,\ldots,g$, $\zik=1$  if row $\ii$ belongs to row-group $\kk$, $0$ otherwise.

In the same way, the latent column assignments $\bw=(\bw_1,\ldots,\bw_d)$ are i.i.d. multinomial variables with $m$ categories: \[\bw_j\sim \mathcal{M}(1, \brho=(\rho_1,\ldots,\rho_m))\] 
For $\ell=1,\ldots,m$, $\wjl=1$  if column $\jj$ belongs to column-group $\el$ and $0$ otherwise.
\item the row  and column assignments  are independent: $p(\bz, \bw)=p(\bz)p(\bw)$

\item conditionally to row and column assignments $\bz\times\bw$, the observed data $\Xij$ are independent, and their conditional  distribution $\dens(.,\alpha)$ belongs to the same parametric family, which parameter $\alpha$ only depends on the given block:
  \begin{equation*}
    \label{eq:distribution-conditionnelle}
    \Xij |\{ \zik\wjl = 1 \}\sim \dens(., \al_{\kk\el}).
  \end{equation*}
\end{itemize}
Hence, the  complete parameter set is $\btheta=(\bpi,\brho,\bal)\in\bTheta$, with $\bal=(\al_{11}, \dots, \al_{\g\m})$ and  $\bTheta$  the parameter space.  Figure \ref{Fig:notations} summarizes these notations. 

\begin{rem}  Group, class and cluster in one hand,    label and assignment in the other hand will be used indistinctly. Moreover, for notation convenience, $\sum_{\ii}$, $\sum_{\jj}$, $\sum_{\kk}$, $\sum_{\ell}$ stand for $\sum_{\ii=1}^n$, $\sum_{\jj=1}^d$, $\sum_{\kk=1}^g$, $\sum_{\ell=1}^m$.
\end{rem}

\begin{figure}[!h]
\input{figure}
\caption{\label{Fig:notations} Notations. Left: Notations for the elements of observed data matrix are in black, notations for the block clusters are in blue. Right: Notations for the model parameter.}
\end{figure}

When performing inference from data, we denote $\bthetas = (\bpis, \brhos, \bals)$ the true parameter set, \emph{i.e.} the parameter values used to generate the data, and $\bzs$ and $\bws$ the true (and usually unobserved) row and column assignments.
 For  indicator membership variables $\bz$ and $\bw$, we also denote:
\begin{itemize}
\item $\zsumk = \sum_{\ii} \zik$ and $\wsuml = \sum_{\jj} \wjl$
\item $\zssumk$ and $\wssuml$ their counterpart for $\bzs$ and $\bws$.
\end{itemize}
  The confusion matrix allows one to compare the partitions.
\begin{dof}[confusion matrices]
  \label{def:confusion}
  For given assignments $\bz$ and $\bzs$ (resp. $\bw$ and $\bws$), we define the \emph{confusion matrix} between $\bz$ and $\bzs$ (resp. $\bw$ and $\bws$), denoted $\Rgbz$ (resp. $\Rmbw$), as follows:
  \begin{equation*}
    \label{eq:confusion-matrices}
    \Rgbz_{\kk\kp} = \frac{1}{\n} \sum_{\ii} \z^\vrai_{\ii\kk} \z_{\ii\kp} \quad \text{and} \quad \Rmbw_{\el\lp} = \frac{1}{\dd} \sum_{\jj} \w^\vrai_{\jj\el} \w_{\jj\lp}
  \end{equation*}
\end{dof}

\subsection{Likelihood}
When the labels are known, the {\em complete log-likelihood} is given by:
\begin{equation}
  \label{eq:log-vraisemblance-complete}
  \begin{aligned}
    \Lc(\bz,\bw;\btheta)&= \log\prob(\bx,\bz,\bw;\btheta)\\
    &= \log\left\{ \left(\prodik{\pik^{\zik}}\right)\left(\prodjl{\rhol^{\wjl}}\right)\left(\prodijkl{\dens\left(\xij;\alkl\right)^{\zik\wjl}}\right) \right\}\\
    &= \log\left\{ \left(\prodi{\pii_{\zi}}\right)\left(\prodj{\rhoo_{\wj}}\right)\left(\prodij{\dens\left(\xij;\al_{\zi\wj}\right)}\right) \right\}.	\\
  \end{aligned}
\end{equation}
In an unsupervised setting,  the labels are  unobserved and the {\em observed log-likelihood} is obtained by marginalization over all the label configurations:
\begin{equation*}
  \label{eq:log-vraisemblance-observee}
  \LL(\btheta)=\log\prob(\bx;\btheta) = \log \left( \sum_{\bz \in \mcZ, \bw \in \mcW}\prob(\bx,\bz,\bw;\btheta) \right).
\end{equation*}
Due to the double missing data structure $\bz$ for rows  and $\bw$ for columns, neither the observed likelihood nor the E-step of the EM algorithm are tractable. Estimation can can nevertheless be performed either by numerical approximation, or by MCMC methods (\cite[see][]{govaert2013co,keribin2015estimation}).

\subsection{Assumptions}\label{sec:assumptions}
We focus here on LBM where $\dens$ belongs to a regular univariate exponential family set in canonical form:
\begin{equation*}
  \dens(x, \al) = b(x)\exp(\al x - \norm(\al)),
\end{equation*}

The canonical parameter $\al$ belongs to a space  $\mathcal{A}$, so that $\dens(\cdot, \al)$ is well defined for all $\al \in \mathcal{A}$. Classical properties of exponential families ensure that $\norm$ is convex, infinitely differentiable on $\mathring{\mathcal{A}}$, and $\normpm$ is well defined on $\normp(\mathring{\mathcal{A}})$. When $X_{\al} \sim \dens(., \al)$, 
\[\Esp[X_\al] = \normp(\al) \hbox{ and } \Var[X_{\al}] = \normp'(\al).\]

Notice that the definition of the exponential family used here relies on an exhaustive statistic that is $X$ itself. This for a simple convenience. Family sets of the form $ \dens(x, \al) = b(x)\exp(\al t(x) - \norm(\al))$ can also be considered, all the further developments as Bernstein  and concentration inequalities then concerning the exhaustive statistics $t(X)$.

%\subsection{Hypothèses sur le paramètre}
Moreover, we make the following assumptions on the parameter space~:
\begin{enumerate}
 % \item $\bTheta$ is compact \textcolor{red}{is it necessary?};
\item[$H_1$]: There exists a positive constant $c$, and a compact $C_\al$ such that 
  \begin{equation*}
%   \label{eq:assumptions}
    \bTheta \subset [c, 1-c]^{\g} \times  [c, 1-c]^{\m} \times C_{\al}^{\g \times \m} \quad \text{with} \quad C_\al \subset \mathring{\mathcal{A}}.
  \end{equation*}
\item[$H_2$]: The true parameter $\bthetas = (\bpis, \brhos, \bals)$ lies in the relative interior of $\bTheta$.
\item [$H_3$]: The mixture  measure of LBM is identifiable: $\bthetas$ is identifiable up to a permutation of the row-labels and column-labels (see definition \ref{def:equivalence} of equivalent parameters).

%\item[$H_3$]: The map $\alpha \mapsto \dens(\cdot, \al)$ is injective.
%\item[$H_4$]: Each row and each column of $\bals$ is unique.
%\item The ratio $\log(\dd)/\n$ (resp. $\log(\n)/\dd$) tend to 0 with $\n$ and $\dd$.
\end{enumerate}
%\begin{rem}

The previous assumptions are standard. Notice that the following conditions are necessary  for $H_3$ to hold:
\begin{itemize}
\item[] $H_{3a}$: The map $\alpha \mapsto \dens(\cdot, \al)$ is injective.
\item[] $H_{3b}$: Each row and each column of $\bals$ is unique.
\end{itemize}
\cite{keribin2015estimation}  gives sufficient conditions for the generic identifiability of the categorical LBM, \textit{i.e.} except on a manifold 
set of null Lebesgue measure in $\bTheta$ and this property is easily extended to the case of observations from a univariate exponential family. For binary SBM, \cite{bickel2013asymptotic} added the  
assumption $p\geq (\log n)/n$ on the parameter $p$ of the Bernoulli distribution  to take into account sparsity.

Assumption~$H_1$ ensures that the group proportions $\pi_{k}$ and $\rho_\ell$ are bounded away from $0$ and $1$ so that no group disappears when $\n$ and $\dd$ go to infinity. It also ensures that $\al$ is bounded away from the boundaries of $\mathcal{A}$ and that there exists a positive value $\kappa>0$, such that $[\al - \neighborsize, \al + \neighborsize] \subset \mathring{\mathcal{A}}$ for all parameters $\al$ of $C_\al$, which is essential to prove a uniform Bernstein inequality on the $X$.
% Assumptions $H_3$ and $H_4$ are necessary to ensure that the model is identifiable. If the map $\alpha \mapsto \dens(., \al)$ is not injective, the model is trivially not identifiable. Similarly, if rows $\kk$ and $\kp$ are identical, we can build a more parsimonious model that induces the same distribution of $\bx$ by merging groups $\kk$ and $\kp$. In the following, we  consider that $g$ and $m$, row- and column- classes (or groups) counts are known.

%\end{rem} 
Moreover, we define the quantity $\delta(\bal)$ that captures the separation
between  row-groups or column-groups: low values of $ \delta(\bal)$ mean that two row-classes or two column-classes are very similar. 
\begin{dof}[class distinctness]
  \label{def:group-distinctness}
  For $\btheta = (\bpi, \brho, \bal) \in \bTheta$. We define:
  \begin{equation*}
    \delta(\bal) = \min\left\{ \min_{\el \neq \lp} \max_{\kk} \KL(\alkl, 
\al_{\kk\lp}), \min_{\kk \neq \kp} \max_{\el} \KL(\alkl, \al_{\kp\el}) \right\}
  \end{equation*}
  with $\KL(\al,\al') = \Esp_{\al}[\log(\dens(X, \al)/\dens(X, \al'))]=\normp(\al) (\al - \al') + \norm(\al') - \norm(\al)$ the Kullback divergence between $\dens(., \al)$ and $\dens(., \al')$.
\end{dof}
\begin{rem}
Since $\bals$ has distinct rows and distinct columns ($H_3$), $\delta(\bals) > 0$.
\end{rem}

\begin{rem}
These assumptions are satisfied for many distributions, including but not limited to:
\begin{itemize}
\item Bernoulli, when the proportion $p$ is bounded away from $0$ and $1$, or natural parameter $\al = \log(p / (1 - p))$ bounded away from $\pm \infty$;
\item Poisson, when the mean $\lambda$ is bounded away from $0$ and $+\infty$, or natural parameter $\al = \log(\lambda)$ bounded away from $\pm \infty$;
\item Gaussian with known variance when the mean $\mu$, which is also the natural parameter, is bounded away from $\pm \infty$.
\end{itemize}
In particular, the conditions stating that $\norm$ is twice differentiable and 
that $\normpm$ exists are equivalent to assuming that $\X_{\al}$ has 
positive and finite variance for all values of $\al$ in the parameter 
space.
\end{rem}

%\begin{rem}
%For the sake of simplicity, we only consider data matrices where most elements $\Xij$ are non-zero. Similar results are available for "zero-inflated" distributions $\X_{\al} \sim \kappa_{\n, \dd}\delta_0 + (1 - \kappa_{\n, \dd})\dens(., \al)$ where  $\kappa_{\n, \dd} \xrightarrow[\n,\dd \to +\infty]{} 0$. If missing values are treated as $0$, the two configurations correspond to dense and sparse data tables. Results concerning zero-inflated distributions are detailed in Section~\ref{sec:sparsity}.
%\end{rem}

\subsection{Model Symmetry}
\label{sec:definitions}
The LBM is a generalized mixture model and as such is subject to label 
switching. 
Moreover, the study of the asymptotics  will involve the complete likelihood where 
%lead to take into account 
symmetry properties on the parameter must be taken into account. We first recall the definition of a permutation in LBM, then define equivalence relationships for assignments and parameter, and discuss model symmetry. 

\begin{dof}[permutation]
  \label{def:permutation}
  Let $s$ be a permutation on $\{1,\dots,\g\}$ and $t$ a permutation on $\{1,\dots,\m\}$. If $\boldsymbol{A}$ is a matrix with $\g$ columns, we define $\boldsymbol{A}^s$ as the matrix obtained by permuting the columns of $\boldsymbol{A}$ according to $s$, \emph{i.e.} for any row $\ii$ and column $\kk$ of $\boldsymbol{A}$, ${A}^s_{\ii \kk} = A_{\ii s(\kk)}$. If $\boldsymbol{B}$ is a matrix with $\m$ columns and $\boldsymbol{C}$ is a matrix with $\g$ rows and $\m$ columns, $\boldsymbol{B}^t$ and $\boldsymbol{C}^{s,t}$ are defined similarly:
  \begin{equation*}
    \boldsymbol{A}^s = \left( A_{\ii s^{}(\kk)} \right)_{\ii,\kk} \quad  \boldsymbol{B}^t = \left( B_{\jj t^{}(\el)} \right)_{\jj, \el} \quad \boldsymbol{C}^{s,t} = \left( C_{s^{}(\kk) t^{}(\el)} \right)_{\kk,\el}
  \end{equation*}
\end{dof}

\begin{dof}[equivalence]
  \label{def:equivalence}
  We define the following equivalence relationships:
  \begin{itemize}
  \item Two assignments $(\bz, \bw)$ and $(\bz', \bw')$ are \emph{equivalent}, denoted $\sim$, if they are equal up to label permutation, \emph{i.e.} there exist two permutations $s$ and $t$ such that $\bz' = \bz^s$ and  $\bw' = \bw^t$.
  \item Two parameters $\btheta$ and $\btheta'$ are \emph{equivalent}, denoted $\sim$, if they are equal up to label permutation, \emph{i.e.} there exist two permutations $s$ and $t$ such that $(\bpi^s, \brho^t, \bal^{s,t}) = (\bpi', \brho', \bal')$. This is \emph{label-switching}.
  \item  $(\btheta, \bz, \bw)$ and $(\btheta', \bz', \bw')$ are \emph{equivalent}, denoted $\sim$, if they are equal up to label permutation on $\bal$, \emph{i.e.} there exist two permutations, $s$ and $t$ such that $(\bal^{s,t}, \bz^s, \bw^{t}) = (\bal', \bz', \bw')$.
  \end{itemize}
\end{dof}

The last equivalence relationship is not concerned with $\bpi$ and $\brho$. It is useful when dealing with the conditional likelihood $\prob(\bx| \bz, \bw; \btheta)$ which depends neither on $\bpi$ nor $\brho$:  in fact, if $(\btheta, \bz, \bw) \sim (\btheta', \bz', \bw')$, then for all $\bx$, we have $\prob(\bx| \bz, \bw; \btheta) = \prob(\bx| \bz', \bw'; \btheta')$. 
Note also that $\bz \sim \bzs$ (resp. $\bw \sim \bws$) if and only if   there exists a permutation of the rows of the confusion matrix $\Rgbz$ (resp. $\Rmbw$)  leading to a diagonal matrix. 

\begin{dof}[symmetry]
  \label{def:symmetry}
  We say that the parameter $\btheta$ \emph{exhibits symmetry for the permutations} $s,t$ if
  \begin{equation*}
    (\bpi^s, \brho^t, \bal^{s,t}) = (\bpi, \brho, \bal).
  \end{equation*}
  $\btheta$ \emph{exhibits symmetry} if it exhibits symmetry for any non trivial pair of permutations $(s,t)$. Finally the set of pairs $(s,t)$ for which $\btheta$ exhibits symmetry is denoted $\Symmetric(\btheta)$.
\end{dof}

%\begin{dof}[$\btheta$-symmetry]
%  \label{def:theta-symmetry}
%  For any $\btheta \in \bTheta$, we say that assignments $(\bz, \bw)$ and $(\bz', \bw')$ are \emph{$\btheta$-equivalent}, noted $\thetasim$, if there exists $(s,t) \in \Symmetric(\btheta)$ such that $(\bz^s, \bw^t) = (\bz', \bw')$. \CK{ ou est-ce qu'on s'en sert?}
%\end{dof}

\begin{rem}
  The set of parameters that exhibit symmetry is a manifold of null Lebesgue measure in $\bTheta$.  This notion of symmetry %allows to deal with the  non-identifiability of the class labels  and 
is subtler than and different from label switching.
To emphasize the difference between equivalence and symmetry, consider the 
following model: $\bpi = (1/2, 1/2)$, $\brho = (1/3, 2/3)$ and $\bal = \left( 
\begin{array}{cc} \al_1 & \al_2 \\ \al_2 & \al_1 \end{array} \right) $ with 
$\al_1 \neq \al_2$. The only permutations of interest here are $s = t = [1\ 2]$. 
Choose any $\bz$ and $\bw$. Because of label switching, we know that $\prob(\bx, 
\bz^s, \bw^t; \btheta^{s,t}) = \prob(\bx, \bz, \bw; \btheta)$. $(\bz^s, \bw^t)$ 
and $(\bz, \bw)$ have the same likelihood but under \emph{different} parameters 
$\btheta$ and $\btheta^{s,t}$. If however, $\brho = (1/2, 1/2)$, then $(s, t) 
\in \Symmetric(\btheta)$ and $\btheta^{s,t} = \btheta$ so that $(\bz, \bw)$ and 
$(\bz^s, \bw^t)$ have the same likelihood under the \emph{same} 
parameter $\btheta$. In particular, if $(\bz, \bw)$ is a maximum-likelihood 
assignment under $\btheta$, so is $(\bz^s, \bw^t)$. In other words, if $\btheta$ 
exhibits symmetry,  the maximum-likelihood assignment is \emph{not unique} 
under the true model and there are at least $\# \Symmetric(\btheta)$ of them. 
This has important implications for the asymptotics of the observed 
likelihood ratio. 
\end{rem}

\subsection{Distance and local assignments} We define the distance up to equivalence between two sets of assignments as follows:
\begin{dof}[distance]
  \label{def:equivalence-distance}
 The distance, up to equivalence, between configurations $\bz$ and $\bzs$ is defined as
    \begin{equation*}
    \|\bz - \bzs\|_{0, \sim} = \inf_{\bz' \sim \bz} \|\bz' - \bzs\|_0
    \end{equation*}
 where, for all matrix $\bz$,  $\left\|\cdot\right\|_{0}$ is the Hamming norm 
\[\left\|\bz\right\|_{0}=\sum_{\ii,\kk}{\mathds{1}{\left\{\zik\neq0\right\}}}.\]
A similar definition is set for the distance between $\bw$ and $\bws$.
\end{dof}
This allows us to define a      neighborhood of radius $r$ in the assignment space, taking 
into account equivalent assignments classes.
\begin{dof}[Set of local assignments]
  \label{def:set-local-assignments}
  We denote $S(\bzs, \bws, r)$ the set of configurations that have a representative (for $\sim$) within relative radius $r$ of $(\bzs, \bws)$:
  \begin{equation*}
    S(\bzs, \bws, r) = \left\{ (\bz, \bw) : \|\bz - \bzs\|_{0, \sim} \leq r \n \text{ and } \|\bw - \bws\|_{0, \sim} \leq r \dd \right\}
  \end{equation*}
\end{dof}

%%%%%%%%%%%%%%%%%%%%%%%%%%%%%

%%%%%%%%%%%%%%%%%%%%%%%%%%%%%%%%%%%%%%%%%%
%%%%% Complete model
%%%%%%%%%%%%%%%%%%%%%%%%%%%%%%%%%%%%%%%%%%
\section{Asymptotic properties in the complete data model}
\label{sec:mle-complete-likelihood}
As stated in the introduction, we first study the asymptotic properties of the complete data model.
Let $\widehat{\btheta}_{c}=\left(\hbpi,\hbrho,\hbal\right)$ be the MLE of $\btheta$ in the complete data model, where the real assignments $\bz=\bzs$ and $\bw=\bws$ are known. We can derive the following general estimates from Equation~\eqref{eq:log-vraisemblance-complete}:

\begin{equation}
  \label{eq:mle-complete-likelihood}
  \begin{aligned}
%    \hpi_{\kk}  = 
    \pizk = \frac{\zsumk}{\n} & \quad
    %\hrho_{\el}  = 
    \rhowl = \frac{\wsuml}{\dd} \\
    \hxklzw = \frac{\sum_{\ii\jj} \xij \zik \wjl}{\zsumk \wsuml} & \quad
    \hal_{\kk\el}  = \alzwkl = \normpm \left(  \hxklzw \right)
  \end{aligned}
\end{equation}

\begin{proposition}
  \label{prop:mle-asymptotic-normality}
  The matrices $\Sigma_{\bpis} = \Diag(\bpis) - \bpis\left(\bpis\right)\transpose $, $\Sigma_{\brhos} = \Diag(\brhos) - \brhos\left(\brhos\right)\transpose $ are semi-definite positive, of rank $\g-1$ and $\m-1$, and $\widehat{\bpi}$ and $\widehat{\brho}$ are asymptotically normal:
  \begin{equation}
    \label{eq:mle-proportion-asymptotic-normality}
    \sqrt{\n}\left( \widehat{\bpi}\left(\bzs\right) - \bpis \right) \xrightarrow[\n \to \infty]{\mathcal{D}} \mathcal{N}(0, \Sigma_{\bpis})
    \quad \text{and} \quad
    \sqrt{\dd}\left( \widehat{\brho}\left(\bws\right) - \brhos \right) \xrightarrow[\dd \to \infty]{\mathcal{D}} \mathcal{N}(0, \Sigma_{\brhos})
  \end{equation}
Similarly, let $V(\bals)$ be the matrix defined by $[V(\bals)]_{\kk\el} = 1/\normp'(\al^\vrai_{\kk\el})$ and\\  $\Sigma_{\bals} = \Diag^{-1}(\bpis) V(\bals) \Diag^{-1}(\brhos)$. Then:
  \begin{equation*}
    \label{eq:mle-parameters-asymptotic-normality}
    \sqrt{\n\dd}\;(\hal_{\kk\el}\left(\bzs,\bws\right) - \al^\vrai_{\kk\el}) \xrightarrow[\n, \dd \to \infty]{\mathcal{D}}  \mathcal{N}(0, \Sigma_{\bals, \kk\el})\;\;\; \hbox{for all  } \kk, \el
  \end{equation*}
  and the components $\hal_{\kk,\el}$ are independent. 
\end{proposition}

\textit{Proof:}
Since $\widehat{\bpi}\left(\bzs\right) = \left(\hpi_1\left(\bzs\right), \dots, \hpi_\g\left(\bzs\right)\right)$ (resp. $\widehat{\brho}\left(\bws\right)$) is the sample mean of $\n$ (resp. $\dd$) i.i.d. multinomial random variables with parameters $1$ and $\bpis$ (resp. $\brhos$), a simple application of the central limit theorem (CLT) gives:
\begin{equation*}
  \Sigma_{\bpis, \kk \kp} =
  \begin{cases}
    \pii^{\vrai}_{\kk}(1 - \pii^{\vrai}_{\kk}) & \text{if} \quad \kk = \kp \\
    -\pii^{\vrai}_{\kk} \pii^{\vrai}_{\kp} & \text{if} \quad \kk \neq \kp \\
  \end{cases}
  \quad \text{and} \quad
  \Sigma_{\brhos, \el \lp} =
  \begin{cases}
    \rhoo^{\vrai}_{\el}(1 - \rhoo^{\vrai}_{\el}) & \text{if} \quad \el = \lp \\
    -\rhoo^{\vrai}_{\el} \rhoo^{\vrai}_{\lp} & \text{if} \quad \el \neq \lp \\
  \end{cases}
\end{equation*}
which proves Equation~\eqref{eq:mle-proportion-asymptotic-normality}
where $\Sigma_{\bpis}$ and $\Sigma_{\brhos}$ are
semi-definite positive of rank $\g - 1$ and $\m - 1$.

Similarly, $\normp\left(\hal_{\kk\el}\left(\bzs,\bws\right)\right)$ is
the average of $\zsumk^{\vrai}\wsuml^{\vrai} =
\n\dd\hpi_\kk\left(\bzs\right)\hrho_\el\left(\bws\right)$
i.i.d. random variables with mean
$\normp\left(\al^\vrai_{\kk\el}\right)$ and variance
$\normp'\left(\al^\vrai_{\kk\el}\right)$. $\n\dd\hpi_\kk\left(\bzs\right)\hrho_\el\left(\bws\right)$
is itself random but
\mbox{$\hpi_\kk\left(\bzs\right)\hrho_\el\left(\bws\right)
  \xrightarrow[\n,\dd \to +\infty]{} \pii^\vrai_{\kk}
  \rhoo^\vrai_{\el}$} almost surely. Therefore, by Slutsky's lemma and
the CLT for random sums of random variables \cite{Shanthikumar1984}, we have:
\begin{align*}
\sqrt{\n\dd\pii^\star_{\kk}\rhoo^\star_{\el}} \left(\normp\left(\hal_{\kk\el}\left(\bzs,\bws\right)\right) - \normp(\al^\vrai_{\kk\el}) \right)  &= \sqrt{\n\dd\pii^\star_{\kk}\rhoo^\star_{\el}} \left( \frac{\sum_{\ii\jj} \Xij \zik^{\vrai} \wjl^{\vrai}}{nd\hpi_\kk\left(\bzs\right)\hrho_\el\left(\bws\right)} - \normp(\al^\vrai_{\kk\el}) \right) \\
& \xrightarrow[\n,\dd \to +\infty]{\mathcal{D}} \mcN\left(0, \normp'(\al^\vrai_{\kk\el})\right)
\end{align*}
The differentiability of $\normpm$ and the delta method then gives:
\begin{equation*}
  \sqrt{\n\dd} \left(\hal_{\kk\el}\left(\bzs,\bws\right) - \al^\vrai_{\kk\el} \right) \xrightarrow[\n,\dd \to +\infty]{\mathcal{D}} \mcN\left(0, \frac{1}{\pii^\star_{\kk}\rhoo^\star_{\el} \normp'(\al^\vrai_{\kk\el})}\right)
\end{equation*}
and the independence results from the independence of $\hal_{\kk\el}\left(\bzs,\bws\right)$ and
$\hal_{\kp\lp}\left(\bzs,\bws\right)$ as soon as $\kk \neq \kp$ or $\el
\neq \lp$, as they involve different sets of independent variables.

\begin{flushright}
$\square$
\end{flushright}

Moreover, the complete model is locally asymptotically normal (LAN), as stated in the following proposition. Note that the unusual condition for $s$ and $t$ arises from the constraints $\bpi^\transpose \mathbf{1}_g = \brho^\transpose \mathbf{1}_m = 1$, where $\mathbf{1}_g$ is the vector of size $g$ filled with $1$, which must be satisfied even after perturbing $\bpi$ (resp. $\brho$ ) with $s$ (resp. $t$).

% \CK{Moreover, the complete model is locally asymptotically normal (LAN), as stated in the following proposition. To state it, we need to get rid of the constraint on $\bpi$ and $\brho$, for example by replacing $\pi_g$ and $\rho_m$ by $1-\sum_{k=1}^{g-1} \pi_k$ and $1-\sum_{\ell=1}^{m-1} \rhol$ respectively. We then define $\tilde\bpi=(\pi_1,\ldots,\pi_{g-1})$, $\tilde\brho=(\rho_1,\ldots,\rho_{m-1})$, $\tilde\Sigma_{\bpis}=\Sigma_{\bpis}[1:g-1;1:g-1]$ and $\tilde\Sigma_{\brhos}=\Sigma_{\brhos}[1:m-1;1:m-1]$ where $\Sig_{\bpis}$, $\Sig_{\brhos}$ and $\Sig_{\bals}$ are defined in Proposition~\ref{prop:mle-asymptotic-normality}. We also define $\tilde\Sig_{\bals}$ %=\Diag(\Sig_{\bals}) as the diagonal $(gm \times gm)$ matrix which diagonal is formed by the coefficients of $\Sig_{\bals}$.}

\begin{proposition}[Local asymptotic normality]\label{prop:LocalAsymp} 
Let $\Lcs$ the map defined by $\btheta = \left(\bpi,\brho,\bal\right)\mapsto \log\prob\left(\bx,\bzs,\bws;\btheta\right)$ and note $I_{\bpis} = \Diag^{-1}(\bpis)$, $I_{\brhos} = \Diag^{-1}(\brhos)$ and $I_{\bals}$ the component-wise inverse of $\Sig_{\bals}$. For any $s$, $t$ and $u$ in a compact set , such that $t^\transpose \mathbf{1}_g = 0$ and $s^\transpose \mathbf{1}_m = 0$, we have:
\begin{align*}
\Lcs\left(\bpis+\frac{s}{\sqrt{\n}},\brhos+\frac{t}{\sqrt{\dd}},\bals+\frac{u}{\sqrt{\n\dd}}\right)
&=\Lcs\left(\bthetas\right) + s\transpose\bY_{\bpis}+t\transpose\bY_{\brhos}+ \Trace(u \transpose \bY_{\bals})\\
& -\left(\frac{1}{2}s\transpose I_{\bpis}s +\frac{1}{2}t\transpose I_{\brhos}t+\frac{1}{2} \Trace[(u \odot u)\transpose I_{\bals}] \right)\\
&+o_{P}(1)\\
\end{align*}
%where $\odot$ denotes the Hadamard product of two matrices (element-wise product) and $\Sig_{\bpis}$, $\Sig_{\brhos}$ and $\Sig_{\bals}$ are defined in Proposition~\ref{prop:mle-asymptotic-normality}.
%
where $\odot$ denotes the Hadamard product of two matrices (element-wise product), $\bY_{\bpis}$, $\bY_{\brhos}$ are asymptotically centered Gaussian vectors of sizes $g$ and $m$ with respective variance matrices $I(\bpis)$ and $I(\brhos)$ and $\bY_{\bals}$ is a random matrix of size $g \times m$ with independent Gaussian components $Y_{\bals, kl} \sim {\cal N}(0, I_{\bals,kl})$.
% \begin{proposition}[Local asymptotic normality]\label{prop:LocalAsymp} 
% Let $\Lcs$ the map defined on \CK{$\tilde\bTheta$ by $\tilde\btheta = \left(\tilde\bpi,\tilde\brho,\bal\right)\mapsto \log\prob\left(\bx,\bzs,\bws;\tilde\btheta\right)$}. For any $s$, $t$ and $u$ in a compact set %\CK{of $\R^{g-1}$, $\R^{m-1}$ and $\R^{(g-1)(m-1)}$ respectively}, 
% we have:
% \begin{align*}
% \Lcs\left(\tilde\bpis+\frac{s}{\sqrt{\n}},\tilde\brhos+\frac{t}{\sqrt{\dd}},\bals+\frac{u}{\sqrt{\n\dd}}\right)
% &=\Lcs\left(\tilde\bthetas\right) + s\transpose\bY_{\bpis}+t\transpose\bY_{\brhos}+ u\transpose\bY_{\bals}\\
% & -\left(\frac{1}{2}s\transpose\CK{\tilde\Sig_{\bpis}^{-1}}s+\frac{1}{2}t\transpose\CK{\tilde\Sig_{\brhos}^{-1}}t+\frac{1}{2} \CK{u\transpose \tilde\Sig_{\bals} u}\right)\\
% &+o_{P}(1)\\
% \end{align*}
% %where $\odot$ denotes the Hadamard product of two matrices (element-wise product) and $\Sig_{\bpis}$, $\Sig_{\brhos}$ and $\Sig_{\bals}$ are defined in Proposition~\ref{prop:mle-asymptotic-normality}.
% %
% where $\bY_{\bpis}$, $\bY_{\brhos}$ are asymptotically $g-1$ and $m-1$ Gaussian vectors with zero mean and respective variance matrices \CK{$\tilde\Sig_{\bpis}^{-1}$}, \CK{$\tilde\Sig_{\brhos}^{-1}$} and $\bY_{\bals}$ is a %matrix
% vector of %asymptotically 
% $gm$ independent Gaussian components with zero mean and variance matrix \CK{$\tilde\Sig_{\bals}$}.
\end{proposition}
\proofbegin
By Taylor expansion, and with the condition $s^\transpose \mathbf{1}_g = t^\transpose \mathbf{1}_m = 0$
\begin{eqnarray*}
& &\!\!\!\!\!\!\!\!\!\!\Lcs\left(\bpis+\frac{s}{\sqrt{\n}},\brhos+\frac{t}{\sqrt{\dd}},\bals+\frac{u}{\sqrt{\n\dd}}\right)\\
&=&\Lcs\left(\bthetas\right)+\frac{1}{\sqrt{\n}}s\transpose\nabla{\Lcs}_{\bpi}\left(\bthetas\right) +\frac{1}{\sqrt{\dd}}t\transpose\nabla{\Lcs}_{\brho}\left(\bthetas\right)+\frac{1}{\sqrt{\n\dd}}\text{Tr}\left(u\transpose\nabla{\Lcs}_{\bal}\left(\bthetas\right)\right)\\
&&\quad+\frac{1}{\n}s\transpose\bH_{\bpi}\left(\bthetas\right)s+\frac{1}{\dd}t\transpose\bH_{\brho}\left(\bthetas\right)t +\frac{1}{\n\dd}\text{Tr}\left((u \odot u)\transpose\bH_{\bal}\left(\bthetas\right)\right)+o_{P}(1)\\
\end{eqnarray*}
where $\nabla{\Lcs}_{\bpi}\left(\bthetas\right)$, $\nabla{\Lcs}_{\brho}\left(\bthetas\right)$ and $\nabla{\Lcs}_{\bal}\left(\bthetas\right)$ denote the respective components of the gradient of $\Lcs$ evaluated at $\bthetas$ and $\bH_{\bpi}$, $\bH_{\brho}$ and $\bH_{\bal}$ denotes the conditional hessian of $\Lcs$ evaluated at $\bthetas$. By inspection, $\bH_{\bpi}/\n$, $\bH_{\brho}/\dd$ and $\bH_{\bal}/\n\dd$ converge in probability to constant matrices and the random vectors $\nabla{\Lcs}_{\bpi}\left(\bthetas\right)/\sqrt{\n}$, $\nabla{\Lcs}_{\brho}\left(\bthetas\right)/\sqrt{\dd}$ and $\nabla{\Lcs}_{\bal}\left(\bthetas\right)/\sqrt{\n\dd}$ converge in distribution to Gaussian vectors by the central limit theorem. %This concludes the proof.
\proofend

%%%%%%%%%%%%%%%%%%%%%%%%%%%%%%%%%%
%%%%%%%%%%%%%%  profil likelihood
%%%%%%%%%%%%%%%%%%%%%%%%%%%%%%%%%%
\section{Profile Likelihood}
\label{sec:profile-likelihood}
%%%%%%%%%%%%%%%%%%%%%%%%%%%%%
Our main result compares the observed likelihood ratio $\prob(\bx; \btheta)/\prob(\bx; \bthetas)$ with the complete likelihood  $\prob(\bx, \bzs, \bws; \btheta)/\prob(\bx, \bzs, \bws; \bthetas)$. 
To study the behavior of these likelihoods, we shall  work conditionally to the true configurations $(\bzs,\bws)$ that have enough observations in each row or column group. We therefore define in section \ref{sec:regular} so called \emph{regular} configurations and prove that they occur with high probability. We then introduce in section \ref{sec:cond-and-prof-likelihood} conditional and profile log-likelihood ratios and state some of their properties. 
\subsection{Regular assignments} \label{sec:regular}
\begin{dof}[$c$-regular assignments]
  \label{def:regular}
  Let $\bz \in \mcZ$ and $\bw \in \mcW$. For any $c > 0$, we say that $\bz$ and $\bw$ are c-\emph{regular} if
  \begin{equation*}
    \label{eq:regular-configuration}
    \min_{\kk} \zsumk \geq {c\n} \quad \text{and} \quad \min_{\el} \wsuml \geq {c\dd}.
  \end{equation*}
\end{dof}

In regular configurations, each row-group for example has $\Om(\n)$ members, where $u_\n=\Om(\n)$ if there exists two constant $a, b>0$ such that for $\n$ enough large $a\n \leq u_\n \leq b\n$. $c/2$-regular assignments, with $c$ defined in Assumption $H_1$, have high  $\Prob_{\bthetas}$-probability in the space of all assignments, uniformly over all $\bthetas \in \bTheta$, as stated in Proposition \ref{cor:prob-regular-configurations-star}.

\begin{proposition}
  \label{cor:prob-regular-configurations-star}
  Define $\mcZ_1$ and $\mcW_1$ as the subsets of $\mcZ$ and $\mcW$
  made of $c/2$-regular assignments, with $c$ defined in assumption $H_1$. Denote $\Om_1$ the event $\{ (\bzs, \bws)
  \in \mcZ_1 \times \mcW_1 \}$, then:
  \begin{equation*}
    \Prob_{\bthetas}\left( \bar{\Om}_1 \right) \leq \g \exp\left( -\frac{\n c^2}{2}\right)+ \m \exp\left( -\frac{\dd c^2}{2}\right).
  \end{equation*}
\end{proposition}

Each $\zsumk$ is a sum of $\n$ i.i.d Bernoulli random variables with parameter $\pik \geq \pii_{\min} \geq c$. The proof is straightforward and stems from a simple Hoeffding bound
\begin{equation*}
  \Prob_{\bthetas}\left( \zsumk \leq \n \frac{c}{2} \right)
  \leq
  \Prob_{\bthetas}\left( \zsumk \leq \n \frac{\pik}{2} \right)
  \leq
  \exp\left( - 2\n\left(\frac{\pik}{2}\right)^2 \right)
  \leq
  \exp\left( - \frac{\n c^2}{2} \right)
\end{equation*}
and a union bound over $\g$ values of $\kk$, with similar approach for $\wsuml$.

\subsection{Conditional and profile log-likelihoods}
\label{sec:cond-and-prof-likelihood} 
%We now consider the likelihood $\prob(\bx,\bz,\bw;\btheta)$ when $\bz$ and $\bw$ are unknown. Instead of maximizing this quantity for the correct assignments $(\bzs, \bws)$, we will maximize it in $\btheta$ it for \enquote{well-behaved} test assignments $(\bz, \bw)$ and show that it is a $\smallO_P$ of $\prob(\bx,\bzs,\bws;\bthetas)$ whenever $(\bz, \bw) \nsim (\bzs, \bws)$.
%
%We first introduce few notations.
 Introducing the  conditional log-likelihood ratio
\begin{equation*}\label{eq:fnd} 
  \Fnd(\btheta, \bz, \bw)  = \log \frac{\prob(\bx| \bz,\bw;\btheta)}{\prob(\bx| \bzs,\bws;\bthetas)},
  \end{equation*}
the complete likelihood can be written as follows
\begin{equation*} \label{eqn:fromcomplete-to-fnd}
\prob(\bx, \bz, \bw; \btheta)= \prob(\bz, \bw; \btheta) \prob(\bx| \bzs, \bws; 
\bthetas) \exp(\Fnd(\btheta, \bz, \bw) ).
\end{equation*}
The study of $\Fnd$ will be of crucial importance, as well as its maximum over $\bTheta$. After some definitions, we examine  some useful properties.
%(\ref{prop:profile-likelihood} and \ref{prop:maximum-conditional-likelihood}) and give an upper bound of the expected profile log-likelihood ratio \ref{prop:profile-likelihood-derivative}.

\begin{dof}
\label{def:conditional-profile-likelihood}
The conditional expectation $G$ of $\Fnd$ is defined as:
\begin{equation*}
  \label{eq:conditional-likelihood}
  \begin{aligned}
    \G(\btheta, \bz, \bw) & = \Esp_{\bthetas} \left[ \left. \log \frac{\prob(\bx| \bz,\bw;\btheta)}{\prob(\bx| \bzs,\bws;\bthetas)} \right| \bzs, \bws  \right] =
    \Esp_{\bthetas} \left[ \left. \Fnd(\btheta, \bz, \bw) \right| \bzs, \bws  \right]
  \end{aligned}
\end{equation*}
Moreover, the profile log-likelihood ratio $\Lamb$ and its expectation $\Lambtilde$ are defined as:
\begin{equation*}
  \label{eq:profile-likelihood}
  \begin{aligned}
    \Lamb(\bz, \bw) & = \max_{\btheta} \Fnd(\btheta, \bz, \bw)  \\
    \Lambtilde(\bz, \bw) & = \max_{\btheta} \G(\btheta, \bz, \bw).
  \end{aligned}
\end{equation*}
\end{dof}

\begin{rem}
As  $\Fnd$ and $\G$ only depend on $\btheta$ through $\bal$, we will sometimes replace $\btheta$ with $\bal$ in the expressions of $\Fnd$ and $\G$. Replacing $F_{\n,\dd}$ and $G$ by their profiled version $\Lamb$ and $\Lambtilde$ allows us to get rid of the continuous argument $\bal$ of $F_{\n\dd}$ and to rely instead only on discrete contrasts $\Lamb$ and $\Lambtilde$.
\end{rem}

\paragraph{}Now, Proposition \ref{prop:profile-likelihood} characterizes which values of $\bal$ maximize $\Fnd$ and $\G$ to reach $\Lamb$ and $\Lambtilde$. Propositions \ref{prop:maximum-conditional-likelihood} and \ref{prop:profile-likelihood-derivative} in turn describes properties  of $G$ and $\Lambtilde$ relative to $(\bz,\bw)$.
\begin{proposition}[maximum of $\G$ and $\Lambtilde$ in $\theta$]
\label{prop:profile-likelihood} Let $\widehat\theta_c=(\hxklzw,\pizk,\rhowl)$ be the maximum likelihood estimator of the complete model, as defined in Equation \ref{eq:mle-complete-likelihood}.
Conditionally on  $\bzs, \bws$, define
 the following quantities:
\begin{equation}
  \label{eq:profile-likelihood-notations}
  \begin{aligned}
    \Sal & = (\Salkl)_{\kk\el} = \left( \normp(\alskl) \right)_{\kk\el} \\
    \barxkl(\bz, \bw) & = \Esp_{\bthetas}[\hxklzw | \bzs, \bws] = \frac{\left[ \Rgbz\transpose \Sal \Rmbw \right]_{\kk\el}}{\pizk\rhowl}
  \end{aligned}
\end{equation}
with $\barxkl(\bz, \bw)=0$ for $\bz$ and $\bw$ such that $\pizk=0$ or $\rhowl=0$.
 
Then $\Fnd(\btheta, \bz, \bw)$ (resp. $\G(\btheta, \bz, \bw)$) is maximum in $\bal$ for  $\bal = \hbal(\bz, \bw)$ (resp. $\bar{\bal}(\bz, \bw)$) defined by:
\begin{equation*}
  \hal(\bz, \bw)_{\kk\el} = \normpm ( \hxklzw ) \quad \text{and} \quad \bar{\al}(\bz, \bw)_{\kk\el} = \normpm ( \barxkl(\bz, \bw) ).
\end{equation*}
Hence,
\begin{equation*}
  \begin{aligned}
    \Lamb(\bz, \bw) & = \Fnd(\hbal(\bz, \bw), \bz, \bw) % = \sum_{\kk} \sum_{\el} \pizk \rhowl \nu(\hxklzw) - \sum_{\kk} \sum_{\el} \pizsk \rhowsl \nu(\Salkl) \\
    \\
    \Lambtilde(\bz, \bw) & = \G(\bar{\bal}(\bz, \bw), \bz, \bw) % = \sum_{\kk} \sum_{\el} \pizk \rhowl \nu(\barxklzw) - \sum_{\kk} \sum_{\el} \pizsk \rhowsl \nu(\Salkl).
  \end{aligned}
\end{equation*}
\end{proposition}
Note that although $\bar{x}_{\kk\el} = \mathbb{E}_{\bthetas}\left[\left.\widehat{x}_{\kk\el}\right|\bzs,\bws\right]$, in general %$\bar{\al}_{\kk\el} = \normpm(\barxkl) \neq \normpm(\hxkl) = \widehat{\al}_{\kk\el}$
$\bar{\al}_{\kk\el} \neq \mathbb{E}_{\bthetas}\left[\left.\widehat{\al}_{\kk\el}\right|\bzs,\bws\right]$ by non linearity of $\normpm$. Nevertheless, since $\normpm$ is Lipschitz over compact subsets of $\normp(\mathring{\mathcal{A}})$, with high probability, $|\bar{\al}_{\kk\el} - \widehat{\al}_{\kk\el}|$ and $|\hxkl - \barxkl|$ are of the same order of magnitude. 

\begin{proposition}[maximum of $\G$ and $\Lambtilde$ in $(\theta,\bz,\bw)$]
  \label{prop:maximum-conditional-likelihood}
  Let $\KL(\al,\al') = \normp(\al) (\al - \al') + \norm(\al') - \norm(\al)$ be the Kullback divergence between $\dens(., \al)$ and $\dens(., \al')$ then:
  \begin{equation}
    \label{eq:conditional-likelihood-second-form}
    \G(\btheta, \bz, \bw) = - \n\dd \sum_{\kk,\kp} \sum_{\el,\lp} \Rgbz_{\kk,\kp} \Rmbw_{\el, \lp} \KL(\al^\vrai_{\kk\el}, \al_{\kp\lp}) \leq 0.
  \end{equation}
  Conditionally on the set $\Om_1$ of regular assignments and for $n, d > 2/c$,
\begin{itemize}
\item[(i)] $\G$ is maximized at $(\bals, \bzs, \bws)$ and its equivalence class.
\item[(ii)] $\Lambtilde$ is maximized at $(\bzs, \bws)$ and its equivalence class and $\Lambtilde(\bzs, \bws)=0$.
%\item[(iii)] \CK{The maximum of $\Lambtilde$ (and hence the maximum of $\G$) is well separated.    A ENLEVER}
\end{itemize}

\end{proposition}

Moreover, the maximum of $\Lambtilde$ in $(\bzs,\bws)$  is  well separated, in the sense that there exists a positive gap between $\Lambtilde (\bzs,\bws)$ and any other $\Lambtilde (\bz,\bw)$ for $(\bz,\bw)$ in a close neighborhood of $(\bzs,\bws)$, as stated in the following proposition:
%Property $(iii)$ of Proposition \ref{prop:maximum-conditional-likelihood} is a direct consequence of the local upperbound for $\Lambtilde$ as stated as follows:

\begin{proposition}[Separability for  $\Lambtilde$]
\label{prop:profile-likelihood-derivative}
Conditionally upon $\Om_1$, there exists a positive constant $C$ such that for all $(\bz, \bw) \in S(\bzs, \bws, C)$:
\begin{equation}
 \label{eq:conditional-likelihood-separability}
  \Lambtilde(\bz, \bw) \leq -\frac{c\delta(\bals)}{4} \left( \dd \|\bz - \bzs\|_{0, \sim} + \n \|\bw - \bws\|_{0,\sim} \right)
\end{equation}
Moreover, there exists a positive constant $B(C)$ such that for all $(\bz, \bw) \notin S(\bzs, \bws, C)$ 
\begin{equation}
 \label{eq:conditional-likelihood-upperbound}
  \Lambtilde(\bz, \bw) \leq -  B(C)\; \n\dd
\end{equation}
\end{proposition}

The proofs of these propositions are reported in Appendix \ref{sec:proofPL}. Proof of Proposition~\ref{prop:profile-likelihood} follows from a straightforward calculation, proof of Proposition~\ref{prop:maximum-conditional-likelihood} uses the technical Lemma~\ref{lemme:casdegalite} to characterize the maximum of $G$ and proof of Proposition~\ref{prop:profile-likelihood-derivative} uses regularity properties of the gradient of $\Lambtilde$ to control its behavior near its maximum.

\section{Main Result}
\label{sec:big-theorem}
%%%%%%%%%%%%%%%%%%%%%%
Our main result  matches the asymptotics of complete and observed likelihoods and is the key to prove the consistency of maximum likelihood and variational estimators. It is set under the assumptions described in section \ref{sec:assumptions} and the following asymptotics for the number of rows $n$ and columns $d$:
\begin{equation*} \label{hyp:ratefornd}
  (H_4):\;\;\;\log (\dd)/\n\rightarrow 0 \hbox{ and } 
\log(\n)/\dd\rightarrow 0.
\end{equation*} 

%stated in Theorem \ref{thm:observed-akin-to-complete-general}

\begin{thm}[complete-observed]
  \label{thm:observed-akin-to-complete-general}
  Let $\bx$ be a matrix of $n\times d$ observations of a LBM with true parameter $\bthetas=(\bpis,\brhos,\bals)$ where the number of row-groups $g$ and column-groups $m$ are known, which conditional distribution belongs to a regular univariate exponential family.  The true random and unobserved assignations for rows and columns are denoted $\bzs$ and $\bws$ respectively.
Define $\# \Symmetric(\btheta)$ as the number of pairs of permutations 
$(s,t)$ for which $\btheta$ exhibits symmetry.

If assumptions $H_1$ to $H_4$ are fulfilled, 
then,  the observed likelihood ratio behaves like the 
complete likelihood ratio, up to a bounded multiplicative factor:
  \begin{equation*}
    \frac{\prob(\bx; \btheta)}{\prob(\bx; \bthetas)} = \frac{\# \Symmetric(\btheta)}{\# \Symmetric(\bthetas)} \max_{\btheta' \sim \btheta} \frac{\prob(\bx, \bzs, \bws; \btheta')}{\prob(\bx, \bzs, \bws; \bthetas)}\left(1 + \smallO_P(1)\right) + \smallO_P(1)
  \end{equation*}
  where both $\smallO_P$ are uniform over all $\btheta \in \bTheta$.
\end{thm}
 The maximum over all $\btheta'$ that are equivalent to $\btheta$ stems from the fact that because of label-switching, $\btheta$ is only identifiable up to its $\sim$-equivalence class from the observed likelihood, whereas it is completely identifiable from the complete likelihood as in this latter case, the labels are known. The terms $\#\Symmetric$ are needed to take into account cases where $\btheta$ exhibits symmetry. These were omitted  
 by \cite{bickel2013asymptotic} for SBM,  although they are also needed in this case, see remark \ref{rem:sym-term}. 
When no $\btheta \in \bTheta$ exhibits symmetry, the following corollary is immediately deduced :
\begin{corollaire}
  \label{cor:observed-akin-to-complete-simple-case}
  If $\bTheta$ contains only parameters that do not exhibit symmetry:
  \begin{equation*}
    \frac{\prob(\bx; \btheta)}{\prob\left(\bx; \btheta^{\vrai}\right)} = \max_{\btheta' \sim \btheta} \frac{\prob(\bx, \bzs, \bws; \btheta')}{\prob(\bx, \bzs, \bws; \bthetas)}\left(1 + \smallO_P(1)\right) + \smallO_P(1)
  \end{equation*}
where the $\smallO_P$ is uniform over all $\btheta \in \bTheta$.
\end{corollaire}

%\textcolor{red}{Redondant avec la fin de la section~\ref{sec:cond-and-prof-likelihood}?}
\paragraph{General sketch of the proof.}
The proof relies on the following decomposition of the observed likelihood:
\[
\prob(\bx; \btheta) =
  \sum_{ (\bz, \bw)} \prob(\bx, \bz, \bw; \btheta) =\sum_{ (\bz, \bw)\sim (\bzs, \bws)} \prob(\bx, \bz, \bw; \btheta) +  \sum_{(\bz, \bw)\nsim (\bzs, \bws)} \prob(\bx, \bz, \bw; \btheta).
\]
where the second term shall be proved to be asymptotically negligible. Its control stems from the study of
the conditional log-likelihood $\Fnd$, see Equation \ref{eqn:fromcomplete-to-fnd}.
%
%Hence, the convergence of $\Fnd$ is of crucial importance to study the asymptotic of the observed likelihood. In fact, we shall have to study it for three types of configurations
In fact, the contribution of  configurations that are not equivalent to $(\bzs, \bws)$ leads itself to the study of a global control, and a sharper local control of $\Fnd$. %: we shall have to establish a global control for configurations far from $(\bzs, \bws)$, and a sharper local control for configurations close to $(\bzs, \bws)$.
Hence, the proof relies on the examination of the asymptotic behavior of  $\Fnd$ on three types of configurations that partition $\mcZ\times\mcW$:

\begin{enumerate}
\item  \emph{global control} for  assignations $(\bz, \bw)$  sufficiently far from  $(\bzs, 
\bws)$, \textit{i.e.} such that $\Lambtilde(\bz, \bw)$ 
is of order $\Omega(-\n\dd)$. Proposition \ref{prop:conditional-likelihood-convergence} gives a large 
deviation result for $\Fnd - \Lambtilde(\bz, \bw)$ to prove that $\Fnd$ is also of 
order $-\Omega_P(\n\dd)$. A key point will be the use of Proposition~\ref{prop:concentration-subexponential}, establishing  a specific concentration inequality for sub-exponential variables. In turn, those assignments contribute 
as a $\smallO_{P}(\prob(\bx, \bzs, \bws; \bthetas))$ to the sum 
(Proposition \ref{prop:large-deviations-profile-likelihood}).

\item \emph{local control}: a small deviation result (Proposition~\ref{prop:profile-likelihood-convergence-local}) is needed  to show that the combined contribution of assignments close to but not equivalent to $(\bzs, \bws)$ is also a $\smallO_{P}(\prob(\bx, \bzs, \bws; \bthetas))$ (Proposition \ref{prop:small-deviations-profile-likelihood}). 
\item \emph{equivalent assignments}: Proposition~\ref{prop:equivalent-configurations-profile-likelihood} examines which of the remaining assignments, all equivalent to $(\bzs, \bws)$, contribute to the sum.
\end{enumerate}

Once these propositions proved, the proof is straightforward, as can be seen below.  They are in turn carefully presented and discussed in dedicated subsections as they represents the core arguments and their proofs are themselves postponed to Appendix \ref{sec:proofCVG} for more readability.
%These results are presented in next section \ref{sec:technicalpropositions} and 

 %%%%%%%%%%%%%%%%%%%%%%%%%
\proofbegin
We work conditionally to $\Omega_1$, defined in Proposition 
\ref{cor:prob-regular-configurations-star}, \emph{i.e.}, the high 
probability event that $(\bzs,\bws)$ is a $c/2$-regular assignment. We choose $(\bzs, \bws)
\in \mcZ_1\times \mcW_1$ and a sequence $\tnd$ decreasing to $0$ but satisfying $\tnd^2 \gg \frac{\n+\dd}{ \n\dd}$. This  is possible when $n\rightarrow \infty$ and $d\rightarrow\infty$, and for example with Assumption $(H_4)$. We write:
\begin{align*}
\prob(\bx; \btheta) &=\sum_{ (\bz, \bw)\sim (\bzs, \bws)} \prob(\bx, \bz, \bw; \btheta)\\
& +  \sum_{(\bz, \bw) \notin S(\bzs, \bws, \tnd)} \prob(\bx, \bz, \bw; \btheta)+
\sum_{\substack{(\bz, \bw) \in S(\bzs, \bws, \tnd) \\ (\bz, \bw) \nsim (\bzs, \bws)}} \prob(\bz, \bw, \bx; \btheta)
\end{align*}
According to Proposition~\ref{prop:large-deviations-profile-likelihood}, conditionally to $\Omega_1$ and for $\n,\dd$ large enough that 
$2\sqrt{2\n\dd}\tnd \geq \g\m$, the contribution of  far away assignments is
\begin{equation*}
    \sup_{\btheta \in \bTheta} \sum_{(\bz, \bw) \notin S(\bzs, \bws, \tnd)} \prob(\bz, \bw, \bx; \btheta) = \smallO_P( \prob(\bzs, \bws, \bx; \bthetas) ).
\end{equation*}
Using the separability of $\Lambtilde$ and Assumption $(H_4)$, Proposition~\ref{prop:small-deviations-profile-likelihood} ensures the existence of $C>0$ such that:
\begin{equation*}
    \sup_{\btheta \in \bTheta} \sum_{\substack{(\bz, \bw) \in S(\bzs, \bws, C) \\ (\bz, \bw) \nsim (\bzs, \bws)}} \prob(\bz, \bw, \bx; \btheta) = \smallO_P( \prob(\bzs, \bws, \bx; \bthetas) )
\end{equation*}
Since $\tnd$ decreases to $0$,  Proposition~\ref{prop:small-deviations-profile-likelihood} can be applied for the local configurations belonging to $S(\tnd)$, for $n, d$ large enough. Therefore the observed likelihood ratio reduces to:
\begin{align*}
    \frac{\prob(\bx; \btheta)}{\prob(\bx; \bthetas)} % & = \frac{\displaystyle \sum_{(\bz, \bw) \sim (\bzs, \bws)} \prob(\bx, \bz, \bw; \btheta) + \sum_{(\bz, \bw) \nsim (\bzs, \bws)} \prob(\bx, \bz, \bw; \btheta)}{\displaystyle \sum_{(\bz, \bw) \sim (\bzs, \bws)} \prob(\bx, \bz, \bw; \bthetas) + \sum_{(\bz, \bw) \nsim (\bzs, \bws)} \prob(\bx, \bz, \bw; \bthetas)} \\
    & = \frac{\displaystyle \sum_{(\bz, \bw) \sim (\bzs, \bws)} \prob(\bx, \bz, \bw; \btheta) + \prob(\bx; \bzs, \bws, \bthetas) \smallO_P(1)}{\displaystyle \sum_{(\bz, \bw) \sim (\bzs, \bws)} \prob(\bx, \bz, \bw; \bthetas) + \prob(\bx; \bzs, \bws, \bthetas) \smallO_P(1)} \\
\end{align*}
Proposition~\ref{prop:equivalent-configurations-profile-likelihood} deals with equivalence and symmetry and allows us to conclude
\begin{equation*}
    \frac{\prob(\bx; \btheta)}{\prob(\bx; \bthetas)} = \frac{\# \Symmetric(\btheta)}{\# \Symmetric(\bthetas)} \max_{\btheta' \sim \btheta} \frac{\prob(\bx, \bzs, \bws; \btheta')}{\prob(\bx, \bzs, \bws; \bthetas)}(1 + \smallO_P(1)) + \smallO_P(1).
\end{equation*}

\proofend

%\proofbegin 

%%%%%%%%%%%%%%%%%%%%%%%%%%% début preuve
\begin{rem}\label{rem:sym-term}
As already pointed out, if $\btheta$ exhibits symmetry, the maximum likelihood 
assignment is not unique under $\btheta$, and $\# \Symmetric(\btheta)$ 
terms contribute with the same weight. This was not taken into account by 
\cite{bickel2013asymptotic}, and it is interesting to see why it should be 
also  present for SBM. Recall that SBM has only one set of labels $\bz$. The 
proof relies on the the decomposition
\[
 \prob(\bx; \btheta) = \sum_{\bz} \prob(\bx,\bz;\btheta)= \sum_{\bz'\sim \bzs}  \prob(\bx,\bz';\btheta) + \sum_{\bz'\nsim \bzs}  \prob(\bx,\bz';\btheta)
\]
where the second term of the sum is neglectible  compared to the first term. Now, $\bz'\sim \bzs$ means that there exists a permutation $t: [g]\rightarrow [g]$ such that $\bz'=\bz^t$ and $p(\bx,\bz^t;\btheta)=p(\bx,\bz;\btheta^t)$. The first term is written on Page 1941, Equation (25) in \cite{bickel2013asymptotic} as  
\[
\sum_{\bz'\sim \bzs}  \prob(\bx,\bz';\btheta)=\sum_{\btheta'\sim\btheta} p(\bx,\bzs;\btheta')=(1+o(1))\max_{\btheta'\sim\btheta}p(\bx,\bzs;\theta')
\]
However, the first equality  is not always correct. Actually, we have
\begin{align*}
\sum_{\bz'\sim \bzs}  \prob(\bx,\bz';\btheta)&= \sum_{t: [g]\rightarrow [g]}  \prob(\bx,\bz^{\star, t};\btheta) = \sum_{t: [g]\rightarrow [g]} \prob(\bx,\bzs;\btheta^t)
\end{align*}
Take a special case of symmetry where $\pi=(1/g,\ldots,1/g)$ and $\alpha=(p-q)I_g + q 1_g 1_g^T$. Then we have $\btheta^t=\btheta$ for all $t$. Thus,
\[
\sum_{\bz'\sim \bzs}  \prob(\bx,\bz';\btheta)=g!\;p(\bx,\bzs;\btheta). 
\]
Even for the SBM, we thus have generally:
\[
\sum_{\bz'\sim \bzs}  \prob(\bx,\bz';\btheta)=(1+o(1)) \# \Symmetric(\btheta)\max_{\btheta'\sim\btheta}p(\bx,\bzs;\btheta')
\]
\end{rem}

 %They are gathered  together in sub-section \ref{sec:proofbigth} to achieve the proof of our main result. 

%\subsection{Different asymptotic behaviors}\label{sec:technicalpropositions}

\subsection{Global Control}

% a global deviation bound for $\Fnd$ and shows that  $\left( \Fnd - \G \right)_{+} = \bigO_{P}(\sqrt{\n\dd})$.
A large deviation inequality for configurations $(\bz, \bw)$ far from $(\bzs, \bws)$ is build and used to prove that far away configurations make a small contribution to $\prob(\bx; \btheta)$.
Since we restricted $\al$ in a bounded subset of $\mathring{\mathcal{A}}$, there exists two positive values $M_{\al}$ and $\neighborsize$ such that $C_\al + (-\neighborsize, \neighborsize) \subset [-M_{\al}, M_\al] \subset \mathring{\mathcal{A}}$. Moreover, the variance of $\X_{\al}$ is bounded away from $0$ and $+\infty$:
\begin{equation*}
  \label{eq:condition-variance}
  \sup_{\al \in [-M_{\al}, M_{\al}] } \Var(\X_\al) = \vmax^2 < +\infty \quad \text{and} \quad \inf_{\al \in [-M_{\al}, M_{\al}]} \Var(\X_\al) = \vmin^2 > 0.
\end{equation*}
\begin{proposition}
With the previous notations, if $\al \in C_\al$ and $X_{\al} \sim \dens(., \al)$, then $X_\al$ is sub-exponential with parameters $(\vmax^2, \neighborsize^{-1})$.
\end{proposition}
The latter proposition is a direct consequence of the definition of sub-exponential variables, see Appendix \ref{sec:subexp}.

\begin{proposition}[large deviations of $\Fnd$]
  \label{prop:conditional-likelihood-convergence}
  Let $\Diam(\bTheta) = \sup_{\btheta, \btheta'} \|\btheta - \btheta' 
\|_{\infty}$. For all $\vareps_{\n, \dd} \leq \neighborsize \vmax$ and $\n, \dd$
  \begin{multline}
    \label{eq:conditional-likelihood-convergence}
    \Delta_{\n\dd}^1(\vareps_{\n\dd})\\ = \Prob\left( \sup_{\btheta, \bz, \bw} \left\{ \Fnd(\btheta, \bz, \bw) - \Lambtilde(\bz,\bw) \right\} \geq \vmax \n\dd \Diam(\bTheta) 2\sqrt{2}\vareps_{\n\dd}\left[1 + \frac{\g\m}{2\sqrt{2\n\dd}\vareps_{\n\dd}} \right] \right) \\ 
    \leq \g^{\n} \m^{\dd} \exp\left( - \frac{\n\dd\vareps_{\n\dd}^2} {2}\right)
  \end{multline}
\end{proposition}
In particular, if $\n$ and $\dd$ are large enough that
$2\sqrt{2\n\dd}\vareps_{\n\dd} \geq \g\m$, the previous inequality ensures 
that with high probability, $\Fnd(\btheta, \bz, \bw) - \Lambtilde(\bz,\bw)$ is 
no greater than $\vmax \n\dd \Diam(\bTheta) 4\sqrt{2}\vareps_{\n\dd}$.

The concentration inequality used in \cite{bickel2013asymptotic} to prove an 
analog result for SBM is not sufficient here, as it can be used only for 
upper-bounded observations, which is obviously not the case for all exponential 
families. We instead develop a Bernstein-type inequality for 
sub-exponential variables (Proposition \ref{prop:concentration-subexponential}) to upper bound 
$\Fnd(\btheta, \bz, \bw) - \Lambtilde(\bz,\bw)$. Proposition 
\ref{prop:conditional-likelihood-convergence} relies heavily on this 
Bernstein inequality. A straightforward consequence of this deviation 
bound is that the combined contribution of assignments far away 
 from $(\bzs, \bws)$ to the sum is negligible, 
assuming that the numbers $\n$ of rows  and $\dd$ of columns grow at 
commensurate rates, as stated in the following proposition:
\begin{proposition}[contribution of far away assignments]
  \label{prop:large-deviations-profile-likelihood}
  Assume $\n\to \infty$ and $\dd \to \infty$, 
  and choose $\tnd$ decreasing to $0$ such that $\tnd^2\gg \frac{n+d}{nd}$. Then 
conditionally on $\Omega_1$ and for $\n,\dd$ large enough that 
$2\sqrt{2\n\dd}\tnd \geq \g\m$, we have:
  \begin{equation*}
    \sup_{\btheta \in \bTheta} \sum_{(\bz, \bw) \notin S(\bzs, \bws, \tnd)} \prob(\bz, \bw, \bx; \btheta) =  \prob(\bzs, \bws, \bx; \bthetas)\smallO_P(1)
  \end{equation*}
 where the $\smallO_P$ is uniform in probability over all $\btheta \in\bTheta$.
\end{proposition}

\subsection{Local Control}

Proposition~\ref{prop:conditional-likelihood-convergence} gives deviations of order $\bigO_{P}(\sqrt{\n\dd})$, which are only useful for $(\bz, \bw)$ such that $\G$ and $\Lambtilde$ are large compared to $\sqrt{\n\dd}$. For $(\bz, \bw)$ close to $(\bzs, \bws)$, we need tighter concentration inequalities, of  order $\smallO_P(\n+\dd)$, as follows:

% \begin{proposition}[small deviations $\Fnd$]
%   \label{prop:profile-likelihood-convergence-local}
%   Conditionally upon $\Om_1$, there exists three positive constants $c_1$, $c_2$ 
% and $C$ such that for all $\vareps \leq \neighborsize \vmin^2$, for all $(\bz, 
% \bw) \nsim(\bzs, \bws)$ such that $(\bz, \bw) \in S(\bzs, \bws, C)$:
%   \begin{equation}
%     \label{eq:profile-likelihood-convergence-local}
%     \begin{aligned}
%       \Delta_{\n\dd}^2(\vareps) = \Prob_{\bthetas}\left( \sup_{\btheta} \frac{\Fnd(\btheta, \bz, \bw) - \Lambtilde(\bz, \bw)}{\dd \| \bz - \bzs \|_{0, \sim} + \n \| \bw - \bws \|_{0, \sim}} \geq \vareps \right) \leq \exp \left( - \frac{\n\dd c^2 \vareps^2}{128(c_1 \vmax^2 + c_2 \neighborsize^{-1} \vareps)}\right)
%     \end{aligned}
%   \end{equation}
% \end{proposition}

%
\begin{proposition}[small deviations $\Fnd$]
  \label{prop:profile-likelihood-convergence-local}
  Conditionally upon $\Om_1$, for $\n$ and $\dd$ satisfying $(H_4)$, and for $(\bz,\bw)   \in S(\bzs, \bws, c/4) $ where $c>0$ is  defined in $(H_1)$, we have: 
  \begin{equation*}
    \label{eq:profile-likelihood-convergence-local}
    \sum_{\substack{(\bz, \bw) \in S(\bzs, \bws, c/4) \\ (\bz, \bw) \nsim (\bzs, \bws)}}\frac{\Lamb(\bz, \bw) - 
\Lambtilde(\bzs, \bws)}{\dd \|\bz - \bzs\|_{0, \sim} + \n\|\bw - \bws\|_{0, 
\sim}} = \smallO_P(1)
  \end{equation*}
\end{proposition}

The next proposition uses  Propositions~\ref{prop:maximum-conditional-likelihood} and \ref{prop:profile-likelihood-convergence-local} to show that the combined contribution to the observed likelihood of assignments close to $(\bzs, \bws)$  is also a $\smallO_P$ of  $\prob(\bzs, \bws, \bx; \bthetas)$:

\begin{proposition}[contribution of local assignments]
  \label{prop:small-deviations-profile-likelihood}
  With the previous notations and for $\n$ and $\dd$ satisfying Assumption $(H_4)$, for any $\tilde{c} \leq \min(C, c/4)$ we have:
  \begin{equation*}
    \sup_{\btheta \in \bTheta} \sum_{\substack{(\bz, \bw) \in S(\bzs, \bws, \tilde{c}) \\ (\bz, \bw) \nsim (\bzs, \bws)}} \prob(\bz, \bw, \bx; \btheta) = \smallO_P( \prob(\bzs, \bws, \bx; \bthetas) )
  \end{equation*}
\end{proposition}

\subsection{Equivalent assignments}
It remains to study the contribution of equivalent assignments.
\begin{proposition}[contribution of equivalent assignments]
  \label{prop:equivalent-configurations-profile-likelihood}
  For all $\btheta \in \bTheta$, we have 
  \begin{equation*}
    \sum_{(\bz, \bw) \sim (\bzs, \bws)} \frac{\prob(\bx, \bz, \bw; \btheta)}{\prob(\bx, \bzs, \bws; \bthetas)} = \# \Symmetric(\btheta) \max_{\btheta' \sim \btheta} \frac{\prob(\bx, \bzs, \bws; \btheta')}{\prob(\bx, \bzs, \bws; \bthetas)} (1 + \smallO_P(1))
  \end{equation*}   
  where the $\smallO_P$ is uniform in $\btheta$. 
\end{proposition}
The maximum over $\btheta' \sim \btheta$ accounts for equivalent configurations whereas $\# \Symmetric(\btheta)$ is needed when $\btheta$ exhibits symmetry, as noticed in Remark \ref{rem:sym-term}.

\section{Asymptotics for the Maximum Likelihood (MLE) and Variational (VE) Estimators}
\label{sec:est-asympt-beh}
This section is devoted to the asymptotics of the MLE and VE in the incomplete data model as a consequence of the main result \ref{thm:observed-akin-to-complete-general}.
\subsection{ML estimator}\label{sect:MLE}

%The asymptotic behavior of the maximum likelihood estimator in the incomplete data model is a direct consequence of Theorem~\ref{thm:observed-akin-to-complete-general}.

\begin{thm}[Asymptotic behavior of $\bthetaEMV$]\label{cor:behaviorEMV} Denote $\bthetaEMV$ the maximum likelihood estimator and use the notations of Proposition~\ref{prop:mle-asymptotic-normality}. There exist permutations $s$ of $\{1, \dots, \g\}$ and $t$ of $\{1, \dots, \m\}$ such that
\begin{eqnarray*}
\widehat{\bpi}\left(\bzs\right)-\bpiEMV^{s}=o_{P}\left(\n^{-1/2}\right),&&\widehat{\brho}\left(\bws\right)-\brhoEMV^{t}=o_{P}\left(\dd^{-1/2}\right),\\
\widehat{\bal}\left(\bzs,\bws\right)-\balEMV^{s,t}&=&o_{P}\left(\left(\n\dd\right)^{-1/2}\right).\\
\end{eqnarray*}

%If $\# \Symmetric(\btheta)\neq1$,  $\bthetaEMV$ is still consistent:  there exist permutations $s$ of $\{1, \dots, \g\}$ and $t$ of $\{1, \dots, \m\}$ such that
%\begin{eqnarray*}
%\widehat{\bpi}\left(\bzs\right)-\bpiEMV^{s}=o_{P}\left(1\right),&&\widehat{\brho}\left(\bws\right)-\brhoEMV^{t}=o_{P}\left(1\right),\\
%\widehat{\bal}\left(\bzs,\bws\right)-\balEMV^{s,t}&=&o_{P}\left(1\right). \\
%\end{eqnarray*}
\end{thm}
%Hence, the maximum likelihood estimator for the LBM is consistent and asymptotically normal, with the same behavior as the maximum likelihood estimator in the complete data model when $\theta$ does not exhibit any symmetry.
%The proof is available in section~\ref{annexe:cor:behaviorEMV}.
The proof %in appendix \ref{annexe:cor:behaviorEMV} 
relies on a Taylor expansion of the complete likelihood near its optimum, like in Proposition~\ref{prop:LocalAsymp}, and on our main theorem.

\proofbegin
Note first that unless $\bTheta$ is constrained and with high probability, $\bthetaEMV$ and $\hbtheta_c(\bzs, \bws)$ exhibit no symmetries. Indeed, equalities like $\hxkl 
= \hx_{\kp, \lp}$ have vanishingly small probabilities of being simultaneously true when $\Xij$ is discrete and null when $\Xij$ is continuous.

Note also that $\Lcs$ has a unique maximum at $\hbtheta_c(\bzs, \bws)$. Furthermore, the curvature of $\Lcs$ at $\hbtheta_c(\bzs, \bws)$ with respect to $\bpi$ (resp. $\brho$, $\bal$) converge in probability to $I_{\bpis}$ (resp. $I_{\brhos}$, $I_{\bals}$) defined in Proposition~\ref{prop:LocalAsymp} by consistency of $\widehat{\btheta}_c$. Therefore any estimator $\hbtheta$ bounded away from $\hbtheta_c(\bzs, \bws)$ satisfies $\Lcs(\hbtheta_c(\bzs, \bws)) - \Lcs(\hbtheta) > \Omega_P(1)$. If $\| \hbtheta_c(\bzs, \bws) - \hbtheta \|= o_P(1)$, a Taylor expansion at $\hbtheta_c(\bzs, \bws)$ gives
\begin{align*}
\Lcs(\hbtheta_c(\bzs, \bws)) - \Lcs(\hbtheta) = & \left( \n(\hbpi_c\left(\bzs,\bws\right) - \hbpi)^\transpose I_{\bpis} (\hbpi_c\left(\bzs,\bws\right) - \hbpi) + \right. \\
& \dd(\hbrho_c\left(\bzs,\bws\right) - \hbrho)^\transpose I_{\brhos} (\hbrho_c\left(\bzs,\bws\right) - \hbrho) + \\
& \left. \n\dd \Trace\left[ \left\{(\hbal_c\left(\bzs,\bws\right) - \hbal)\odot(\hbal_c\left(\bzs,\bws\right) - \hbal)\right\}^\transpose I_{\bals} \right] \right) \\
& \times (1 + \smallO_P(1)) + \smallO_P(1)
\end{align*}
where the linear term in the expansion vanishes as $\hbtheta_c(\bzs, \bws)$ is the argmax of $\Lcs$ and the Hessian of $\Lcs$ at $\hbtheta_c(\bzs, \bws)$ were replaced by their limit in probability.

We may now prove the corollary by contradiction. Assume that $\min_{s} (\hbpi_{MLE}^{s} - \widehat{\bpi}\left(\bws\right)) \neq 
\smallO_P\left(\n^{-1/2}\right)$,
$\min_{t} (\hbrho_{MLE}^{t} - \widehat{\brho}\left(\bws\right))\neq 
\smallO_P\left(\dd^{-1/2}\right)$ or 
$\min_{s,t} (\hbal_{MLE}^{s,t} - \widehat{\bal}\left(\bzs,\bws\right)) \neq 
\smallO_P\left(\n\dd^{-1/2}\right)$ where $s$ and $t$ are 
permutations of $\{1, \dots, \g\}$ and $\{1, \dots, \m\}$. Plugging $\hbtheta_{MLE}^{s,t}$ in the previous expansion shows that:
\begin{equation}\label{eq:proof:ConsistenceEMV}
\min_{s, t} \Lcs\left(\widehat{\btheta}_c\left(\bzs,\bws\right)\right)- 
\Lcs\left(\bthetaEMV^{s,t}\right)=\Omega_{P}(1).
\end{equation}
But, since $\widehat{\btheta}_c\left(\bzs,\bws\right)$ and $\bthetaEMV$ maximize 
respectively $\frac{\prob(\bx, \bzs, \bws; \btheta)}{\prob(\bx, \bzs, \bws; 
\bthetas)}$ and $\frac{\prob(\bx; \btheta)}{\prob\left(\bx; 
\btheta^{\vrai}\right)}$ and have no symmetries, it follows by 
Theorem~\ref{thm:observed-akin-to-complete-general} that 
\[\left|\frac{\prob\left(\bx, \bzs, \bws; 
\widehat{\btheta}_c\left(\bzs,\bws\right)\right)}{\prob(\bx, \bzs, \bws; \bthetas)}- 
\max_{s, t} \frac{\prob\left(\bx, \bzs, \bws; 
\bthetaEMV^{s,t}\right)}{\prob(\bx, \bzs, \bws; \bthetas)}\right|=\smallO_P(1)\]
which contradicts Equation~(\ref{eq:proof:ConsistenceEMV}) and concludes the proof.
\proofend
%\CK{ !! et que se passe-t- il si $\theta$ ou $\theta^*$ présente une symmetrie?}
% Si \thetas a une symétrie, on a juste une constante. Sinon, \thetahat ne peut pas avoir de symétries ps.

\subsection{Variational estimator}

Due to the complex dependence structure of the observations, the maximum likelihood estimator of the LBM is not numerically tractable, even  with the EM-algorithm. In practice, a variational approximation can be used, see for example \cite{govaert2003}: for any joint distribution $\setQ\in\mcQ$ on $\mcZ\times\mcW$ a lower bound of $\mcL(\btheta)$ is given by
\begin{eqnarray*}
\Jvar\left(\setQ,\btheta\right)&=&\mcL(\btheta)-KL\left(\setQ,\prob\left(., . ; \btheta,\bx\right)\right)\\
                            &=&\Esp_{\setQ}\left[\Lc\left(\bz,\bw;\btheta\right)\right]+\mcH\left(\setQ\right).
\end{eqnarray*}
where $\mcH\left(\setQ\right)=-\Esp_{\setQ}[\log(\setQ)]$.
Choose $\mcQ$  to be the set of factorized distributions, such that for all $\left(\bz,\bw\right)$
\[\setQ\left(\bz,\bw\right)=\setQ\left(\bz\right)\setQ\left(\bw\right)=\prod_{\ii,\kk}\setQ\left(\zik=1\right)^{\zik} \prod_{\jj,\el}\setQ\left(\wjl=1\right)^{\wjl}\]
allows to obtain tractable expressions of $\Jvar\left(\setQ,\btheta\right)$ as a lower bound of the log-likelihood. The variational estimate $\bthetavar$ of $\btheta$ is defined as
\[\bthetavar\in\underset{\btheta\in\bTheta}{\arg\!\max}\;\underset{\setQ\in\mcQ}{\max}\;\Jvar\left(\setQ,\btheta\right).\]

The following corollary states that $\bthetavar$ has the same asymptotics as $\bthetaEMV$ and $\widehat{\btheta}_{c}$.
%\begin{corollaire}[Variational estimate]\label{cor:Variational}
%Under the conditions of Theorem~\ref{thm:observed-akin-to-complete-general}, \VB{ there exist permutations $s$ of $\{1, \dots, \g\}$ and $t$ of $\{1, \dots, \m\}$ such that
%\begin{eqnarray*}
%\widehat{\bpi}\left(\bzs\right)-\bpivar^{s}=o_{P}\left(\n^{-1/2}\right),&&\widehat{\brho}\left(\bws\right)-\brhovar^{t}=o_{P}\left(\dd^{-1/2}\right),\\
%\widehat{\bal}\left(\bzs,\bws\right)-\balvar^{s,t}&=&o_{P}\left(\left(\n\dd\right)^{-1/2}\right).\\
%\end{eqnarray*}}
%%\[\bthetavar\in
%%\underset{\btheta}{\arg\!\max}\log\prob(\bx;\btheta).\]
%If $\# \Symmetric(\btheta)=1$ and with the notations of Theorem~\ref{thm:observed-akin-to-complete-general}, we have
%  \begin{equation*}
%    \frac{\underset{\setQ\in\mcQ}{\max}\;\exp\left[\Jvar\left(\setQ,\btheta\right)\right]}{\prob(\bx; \bthetas)} =  \frac{1}{\# \Symmetric(\bthetas)} \max_{\btheta' \sim \btheta}\frac{\prob(\bx, \bzs, \bws; \btheta')}{\prob(\bx, \bzs, \bws; \bthetas)}\left(1 + \smallO_P(1)\right) + \smallO_P(1)
%  \end{equation*}
%  where the $\smallO_P$ is uniform over all $\btheta \in \bTheta\backslash\left\{\btheta\in\bTheta|\# \Symmetric(\btheta)=1\right\}$.
%\end{corollaire}
\begin{thm}[Variational estimate]\label{cor:Variational}
Under the assumptions of Theorem~\ref{thm:observed-akin-to-complete-general} there exist permutations $s$ of $\{1, \dots, \g\}$ and $t$ of $\{1, \dots, \m\}$ such that
\begin{eqnarray*}
\widehat{\bpi}\left(\bzs\right)-\bpivar^{s}=o_{P}\left(\n^{-1/2}\right),&&\widehat{\brho}\left(\bws\right)-\brhovar^{t}=o_{P}\left(\dd^{-1/2}\right),\\
\widehat{\bal}\left(\bzs,\bws\right)-\balvar^{s,t}&=&o_{P}\left(\left(\n\dd\right)^{-1/2}\right).
\end{eqnarray*}
%\[\bthetavar\in
%\underset{\btheta}{\arg\!\max}\log\prob(\bx;\btheta).\]
\begin{comment}
Moreover, 
  \begin{eqnarray}  
  \label{cor:Variational:plus}
    \frac{\underset{\setQ\in\mcQ}{\max}\;\exp\left[\Jvar\left(\setQ,\btheta\right)\right]}{\prob(\bx; \bthetas)} =  \frac{1}{\# \Symmetric(\bthetas)} \max_{\btheta' \sim \btheta}\frac{\prob(\bx, \bzs, \bws; \btheta')}{\prob(\bx, \bzs, \bws; \bthetas)}\left(1 + \smallO_P(1)\right) + \smallO_P(1)
\end{eqnarray} 
  where the $\smallO_P$ is uniform over all $\btheta \in \bTheta\backslash\left\{\btheta\in\bTheta|\# \Symmetric(\btheta)=1\right\}$.
  \end{comment}
\end{thm}
\proofbegin
Remark first that for every $\btheta$ and for every $\left(\bz,\bw\right)$, 
\[\prob\left(\bx,\bz,\bw;\btheta\right)\leq\exp\left[\Jvar\left(\delta_{\bz}
\times\delta_{\bw},\btheta\right)\right]
\leq 
\underset{\setQ\in\mcQ}{\max}\;\exp\left[\Jvar\left(\setQ,\btheta\right)\right]
\leq\prob\left(\bx;\btheta\right)\]
where $\delta_{\bz}$ denotes the dirac mass on $\bz$. By dividing by 
$\prob\left(\bx;\bthetas\right)$, we obtain
\[\frac{\prob\left(\bx,\bz,\bw;\btheta\right)}{\prob\left(\bx;\bthetas\right)}
\leq 
\frac{\underset{\setQ\in\mcQ}{\max}\;\exp\left[\Jvar\left(\setQ,
\btheta\right)\right]}{\prob\left(\bx;\bthetas\right)} 
\leq\frac{\prob\left(\bx;\btheta\right)}{\prob\left(\bx;\bthetas\right)}.\]
As this inequality is true for every couple $(\bz,\bw)$, we have in particular:
\[
\underset{(\bz, \bw) \sim 
(\bzs,\bws)}{\max}\frac{\prob\left(\bx,\bz,\bw;\btheta\right)}{ 
\prob\left(\bx;\bthetas\right) } = \underset{\btheta' \sim \btheta}{\max} 
\frac{\prob\left(\bx,\bzs,\bws;\btheta'\right)}{ 
\prob\left(\bx;\bthetas\right) }\leq 
\frac{\underset{\setQ\in\mcQ}{\max}\;\exp\left[\Jvar\left(\setQ, 
\btheta\right)\right]}{\prob\left(\bx;\bthetas\right)} .\]
Noticing that 
$\prob\left(\bx;\bthetas\right) = \#\Symmetric(\bthetas)\prob\left(\bx,
\bzs,\bws;\bthetas\right)(1+o_p(1))$, 
Theorem~\ref{thm:observed-akin-to-complete-general} therefore leads to the 
following bounds:
\begin{multline*}
\underset{\btheta' \sim \btheta}{\max} 
\frac{\prob\left(\bx,\bzs,\bws;\btheta'\right)}{ 
\prob\left(\bx, \bzs, \bws;\bthetas\right)}(1+\smallO_P(1)) \leq 
\frac{\underset{\setQ\in\mcQ}{\max}\;\exp\left[\Jvar\left(\setQ, 
\btheta\right)\right]}{\prob\left(\bx, \bzs, \bws;\bthetas\right)} \\
\leq 
\#\Symmetric(\btheta) \underset{\btheta' \sim \btheta}{\max} 
\frac{\prob\left(\bx,\bzs,\bws;\btheta'\right)}{\prob\left(\bx, \bzs, 
\bws;\bthetas\right)}(1+\smallO_P(1)) + \smallO_P(1).
\end{multline*}
Again, unless $\bTheta$ is constrained, $\bthetaVAR$ exhibits no symmetries with 
high probability and the same proof by contradiction as in 
section~\ref{sect:MLE} gives the result.
\proofend

\section{Conclusion}

The Latent Block Model offers challenging theoretical questions. We solved under mild assumptions the consistency and asymptotic normality of the maximum likelihood and variational estimators for observations with conditional density belonging to a univariate exponential family, and for a balanced asymptotic rate between  the number of rows $n$ and the number of columns $d$:  $\log (\dd)/\n\rightarrow 0$ and $\log(\n)/\dd\rightarrow 0$ as $n$ and $d$ tend to infinity. Our results extend those of \cite{bickel2013asymptotic} for binary SBM not only by managing  the double direction of  LBM, but also  by considering larger types of observations. That brought us to define specific concentration inequalities as large and moderate deviations concerning sub-exponential variables. Moreover, we dealt with specific cases of symmetry that were not taken into account as of now.

A specific framework of sparsity was studied by \cite{bickel2013asymptotic}. This is especially convenient for SBM, to model reasonable network settings: 
increasing sparsity (\textit{i.e.} number of $0$) can be done directly by scaling the   Bernoulli parameters $p_{k\ell}$ with a common factor that should decrease no faster than $\Omega(\log^\delta(n) /n)$, with $\delta > 2$, to ensure consistency. This could also be considered for binary LBM. However this approach  fails to model actual observations in the more general valued setting. The equivalent approach could be to consider the product of a Bernoulli variable with the actual observation value. Note however, than even without considering sparsity we recover essentially the same rate: in the sparse-SBM case, each node should be connected to $\Omega(\log^\delta(n))$ others to ensure consistency whereas in the dense-LBM case, each of the $n$-row should should be characterized by $\Omega(\log^\delta(n))$ columns (and vice-versa) to ensure consistency. 

Alternative research direction could be to  explore asymptotic settings where the numbers $\n$ of rows and $\dd$ columns grow at very different rates. Other open question concern estimation of the number of row and column groups and settings where the number of groups increases with $\n$ and $\dd$.

\appendix

\include{ejs-Theorem1Propositions}

\bibliographystyle{plain}
\bibliography{ejs-BraultKeribinMariadassou}

\end{document}

%% file: commandes.tex
\newcommand{\proofbegin}{\paragraph{Proof.}\hspace{1em}\\}
\newcommand{\proofend}{\begin{flushright}
$\square$
\end{flushright}}

\newcommand{\dd}{d}
\newcommand{\g}{g}
\newcommand{\ii}{i}
\newcommand{\jj}{j}
\newcommand{\kk}{k}
\newcommand{\el}{\ell}
\newcommand{\m}{m}
\newcommand{\n}{n}
\newcommand{\vrai}{\star}
\newcommand{\al}{\alpha}
\newcommand{\als}{\al^\vrai}
\newcommand{\alkl}{\al_{\kk\el}}
\newcommand{\alskl}{\als_{\kk\el}}
\newcommand{\bal}{\boldsymbol{\al}}
\newcommand{\bals}{\boldsymbol{\al}^\vrai}
\newcommand{\etaa}{\eta}

\newcommand{\thetaa}{\theta}
\newcommand{\Thetaa}{\Theta}
\newcommand{\btheta}{\boldsymbol{\thetaa}}
\newcommand{\hbtheta}{\widehat{\btheta}}
\newcommand{\bthetas}{\boldsymbol{\thetaa}^\vrai}
\newcommand{\bTheta}{\boldsymbol{\Thetaa}}
\newcommand{\pii}{\pi}
\newcommand{\bpi}{\boldsymbol{\pii}}
\newcommand{\pik}{\pii_\kk}

\newcommand{\bpis}{\boldsymbol{\pii}^\vrai}

\newcommand{\rhoo}{\rho}
\newcommand{\rhol}{\rhoo_\el}
\newcommand{\brho}{\boldsymbol{\rhoo}}
\newcommand{\brhos}{\boldsymbol{\rhoo}^\vrai}

\newcommand{\w}{w}
\newcommand{\bw}{\mathbf{\w}}
\newcommand{\bws}{\bfw^\vrai}
\newcommand{\wjl}{\w_{\jj\el}}
\newcommand{\x}{x}
\newcommand{\bx}{\mathbf{\x}}
\newcommand{\xij}{\x_{\ii\jj}}
\newcommand{\X}{X}
\newcommand{\Xij}{\X_{\ii\jj}}
\newcommand{\z}{z}
\newcommand{\bfz}{\mathbf{\z}}
\newcommand{\bz}{\bfz}
\newcommand{\bfw}{\bw}

\newcommand{\zik}{\z_{\ii\kk}}
\newcommand{\bzs}{\bfz^\vrai}
\newcommand{\wj}{\w_{\jj}}
\newcommand{\zi}{\z_{\ii}}

\newcommand{\Sal}{\boldsymbol{S}^{\vrai}}

\newcommand{\Salkl}{S^\vrai_{\kk\el}}
\newcommand{\Salmin}{S^\vrai_{\min}}
\newcommand{\Salmax}{S^\vrai_{\max}}

\newcommand{\mcL}{\mathcal{L}}

\newcommand{\Lc}{\mcL_c}
\newcommand{\prodik}{\prod_{\ii,\kk}}
\newcommand{\prodjl}{\prod_{\jj,\el}}
\newcommand{\prodijkl}{\prod_{\ii,\jj,\kk,\el}}
\newcommand{\prodi}{\prod_{\ii}}
\newcommand{\prodj}{\prod_{\jj}}
\newcommand{\prodij}{\prod_{\ii,\jj}}
\newcommand{\prob}{p}
\newcommand{\Prob}{\mathbb{P}}
\newcommand{\Esp}{\mathbb{E}}
\newcommand{\Var}{\mathbb{V}}
\newcommand{\dens}{\varphi}
\newcommand{\norm}{\psi}
\newcommand{\normp}{\norm'}
\newcommand{\normpm}{(\normp)^{-1}}
\newcommand{\mcZ}{\mathcal{Z}}
\newcommand{\mcW}{\mathcal{W}}
\newcommand{\LL}{\mcL}

\newcommand{\Lcs}{\Lc^{\vrai}}

\newcommand{\sommecolonne}{+}
\newcommand{\zsumk}{\z_{\sommecolonne\kk}}
\newcommand{\wsuml}{\w_{\sommecolonne\el}}

\newcommand{\zs}{\z^{\vrai}}
\newcommand{\ws}{\w^{\vrai}}
\newcommand{\zssumk}{\zs_{\sommecolonne\kk}}
\newcommand{\wssuml}{\ws_{\sommecolonne\el}}

\newcommand{\Sig}{\Sigma}

\newcommand{\mcN}{\mathcal{N}}

\newcommand{\hpi}{\widehat{\pii}}
\newcommand{\hbpi}{\widehat{\bpi}}

\newcommand{\hrho}{\widehat{\rhoo}}
\newcommand{\hbrho}{\widehat{\brho}}

\newcommand{\pizk}{\widehat{\pii}_{\kk}(\bz)}
\newcommand{\pizsk}{\widehat{\pi}_{\kk}(\bzs)}
\newcommand{\rhowl}{\widehat{\rhoo}_{\el}(\bw)}
\newcommand{\rhowsl}{\widehat{\rhoo}_{\el}(\bws)}

\newcommand{\hal}{\widehat{\al}}
\newcommand{\hbal}{\widehat{\bal}}

\newcommand{\alzwkl}{\widehat{\al}_{\kk\el}(\bz,\bw)}

\newcommand{\kp}{\kk'}

\def\AArm{\fam0 }%\tenrm}%
\def\R{{\AArm I\!R}}%
\def\un{{\AArm 1\!I}}%

\newcommand{\Rg}{\R_{\g}}
\newcommand{\Rm}{\R_{\m}}
\newcommand{\Rgbz}{\Rg(\bz)}
\newcommand{\Rmbw}{\Rm(\bw)}

\newcommand{\transpose}{{^{T}}}
\newcommand{\vareps}{\varepsilon}

\newcommand{\Lamb}{\Lambda}

\newcommand{\Lambtilde}{\tilde{\Lamb}}
\newcommand{\lp}{\el'}

\newcommand{\barx}{\bar{x}}
\newcommand{\hx}{\widehat{\x}}
\newcommand{\hxkl}{\hx_{\kk\el}}
\newcommand{\hxklzw}{\hx_{\kk\el}(\bz, \bw)}
\newcommand{\baral}{\bar{\al}}
\newcommand{\baralkl}{\baral_{\kk\el}}

\newcommand{\barxkl}{\barx_{\kk\el}}
\newcommand{\barxklzw}{\barx_{\kk\el}(\bz, \bw)}

\newcommand{\neighborsize}{\kappa}
\newcommand{\vmin}{\underline{\sigma}}
\newcommand{\vmax}{\bar{\sigma}}

\newcommand{\tnd}{t_{\n\dd}}

\newcommand{\Fnd}{F_{\n\dd}}
\newcommand{\G}{G}

\newcommand{\bigO}{\mathcal{O}}
\newcommand{\smallO}{o}
\newcommand{\Om}{\Omega}

\DeclareMathOperator{\Diag}{Diag}
\DeclareMathOperator{\Trace}{Tr}
\DeclareMathOperator{\Diam}{Diam}
\DeclareMathOperator{\Symmetric}{Sym}
\DeclareMathOperator{\KL}{KL}

\newcommand{\Y}{Y}
\newcommand{\y}{y}

\newcommand{\bY}{\mathbf{\Y}}
\newcommand{\bH}{\mathbf{H}}
\newcommand{\bthetaEMV}{\widehat{\btheta}_{MLE}}
\newcommand{\bpiEMV}{\widehat{\bpi}_{MLE}}
\newcommand{\brhoEMV}{\widehat{\brho}_{MLE}}
\newcommand{\balEMV}{\widehat{\bal}_{MLE}}
\newcommand{\setQ}{\mathbb{Q}}
\newcommand{\Jvar}{J}
\newcommand{\mcH}{\mathcal{H}}
\newcommand{\bthetavar}{\widehat{\btheta}_{var}}
\newcommand{\balvar}{\widehat{\bal}_{var}}
\newcommand{\brhovar}{\widehat{\brho}_{var}}
\newcommand{\bpivar}{\widehat{\bpi}_{var}}
\newcommand{\mcQ}{\mathcal{Q}}

\newcommand{\bthetaMC}{\widehat{\btheta}_{c}}
\newcommand{\bthetaVAR}{\widehat{\btheta}_{VAR}}

%% file: figure.tex
\begin{center}
\begin{tabular}{cc}
\begin{tikzpicture}[scale=0.85]
%Petit quadrillage
\foreach  \y in {1.5,2,2.5,4.5,6.5}
\draw (\y,4.75) -- (\y,7);
\foreach  \y in {1.5,2,2.5,4.5,6.5}
\draw (\y,1) -- (\y,2.25);
\foreach  \y in {1.5,2,2.5,4.5,6.5}
\draw (\y,2.75) -- (\y,4.25);
\foreach  \y in {1.5,3.5,5.5,6,6.5}
\draw (1,\y) -- (3.25,\y);
\foreach  \y in {1.5,3.5,5.5,6,6.5}
\draw (3.75,\y) -- (5.25,\y);
\foreach  \y in {1.5,3.5,5.5,6,6.5}
\draw (5.75,\y) -- (7,\y);
%Pointillé
\foreach  \y in {1.5,2,2.5,4.5,6.5}
\draw[dotted] (\y,2.25) -- (\y,2.75);
\foreach  \y in {1.5,2,2.5,4.5,6.5}
\draw[dotted] (\y,4.25) -- (\y,4.75);
\foreach  \y in {1.5,3.5,5.5,6,6.5}
\draw[dotted] (3.25,\y) -- (3.75,\y);
\foreach  \y in {1.5,3.5,5.5,6,6.5}
\draw[dotted] (5.25,\y) -- (5.75,\y);

%Bloc bleu
\foreach  \y in {3,4,5,6}
\draw[ultra thick,color=blue] (\y,4.75) -- (\y,7);
\foreach  \y in {3,4,5,6}
\draw[ultra thick,color=blue] (\y,1) -- (\y,2.25);
\foreach  \y in {3,4,5,6}
\draw[ultra thick,color=blue] (\y,2.75) -- (\y,4.25);
\foreach  \y in {2,3,4,5}
\draw[ultra thick,color=blue] (1,\y) -- (3.25,\y);
\foreach  \y in {2,3,4,5}
\draw[ultra thick,color=blue] (3.75,\y) -- (5.25,\y);
\foreach  \y in {2,3,4,5}
\draw[ultra thick,color=blue] (5.75,\y) -- (7,\y);
%Pointillé
\foreach  \y in {3,4,5,6}
\draw[dotted,ultra thick,color=blue] (\y,2.25) -- (\y,2.75);
\foreach  \y in {3,4,5,6}
\draw[dotted,ultra thick,color=blue] (\y,4.25) -- (\y,4.75);
\foreach  \y in {2,3,4,5}
\draw[dotted,ultra thick,color=blue] (3.25,\y) -- (3.75,\y);
\foreach  \y in {2,3,4,5}
\draw[dotted,ultra thick,color=blue] (5.25,\y) -- (5.75,\y);

%Bord
\foreach  \y in {1,7}
\draw[ultra thick,color=black] (\y,4.75) -- (\y,7);
\foreach  \y in {1,7}
\draw[ultra thick,color=black] (\y,1) -- (\y,2.25);
\foreach  \y in {1,7}
\draw[ultra thick,color=black] (\y,2.75) -- (\y,4.25);
\foreach  \y in {1,7}
\draw[ultra thick,color=black] (1,\y) -- (3.25,\y);
\foreach  \y in {1,7}
\draw[ultra thick,color=black] (3.75,\y) -- (5.25,\y);
\foreach  \y in {1,7}
\draw[ultra thick,color=black] (5.75,\y) -- (7,\y);
%Pointillé
\foreach  \y in {1,7}
\draw[dotted,ultra thick,color=black] (\y,2.25) -- (\y,2.75);
\foreach  \y in {1,7}
\draw[dotted,ultra thick,color=black] (\y,4.25) -- (\y,4.75);
\foreach  \y in {1,7}
\draw[dotted,ultra thick,color=black] (3.25,\y) -- (3.75,\y);
\foreach  \y in {1,7}
\draw[dotted,ultra thick,color=black] (5.25,\y) -- (5.75,\y);

% Les indices en ligne
\draw (1,6.75) node[left]{$1$};
\draw (1,6.4) node[left]{$\vdots$};
\draw (1,5.75) node[left]{$\ii$};
\draw (1,5.35) node[left]{$\vdots$};
\draw (1,1.8) node[left]{$\vdots$};
\draw (1,1.25) node[left]{$\n$};
\draw[<->] (0.5,1) -- (0.5,7);
\draw (0.5,4) node[left]{$\n$};
\foreach  \y in {2.3,2.8,3.3,3.8,4.3,4.8}
\draw (1,\y) node[left]{$\vdots$};

% Les indices en colonnes
\draw (1.25,7) node[above]{$1$};
\draw (1.75,7) node[above]{$\cdots$};
\draw (2.25,7) node[above]{$\jj$};
\draw (2.75,7) node[above]{$\cdots$};
\draw (6.25,7) node[above]{$\cdots$};
\draw (6.75,7) node[above]{$\dd$};
\draw[<->] (1,7.5) -- (7,7.5);
\draw (4,7.5) node[above]{$\dd$};
\foreach  \y in {3.25,3.75,4.25,4.75,5.25,5.75}
\draw (\y,7) node[above]{$\cdots$};

%Cases
\draw (1.25,6.5) node[above]{\footnotesize{$\x_{11}$}};
\draw (1.25,5.5) node[above]{\footnotesize{$\x_{\ii1}$}};
\draw (2.25,6.5) node[above]{\footnotesize{$\x_{1\jj}$}};
\draw (2.25,5.5) node[above]{\footnotesize{$\xij$}};
\draw (1.25,1) node[above]{\footnotesize{$\x_{\n1}$}};
\draw (2.25,1) node[above]{\footnotesize{$\x_{\n\jj}$}};
\draw (6.75,6.5) node[above]{\footnotesize{$\x_{1\dd}$}};
\draw (6.75,5.5) node[above]{\footnotesize{$\x_{\ii\dd}$}};
\draw (6.75,1) node[above]{\footnotesize{$\x_{\n\dd}$}};
%Pointillés
\foreach  \yy in {1.25,2.25,6.75}
\foreach  \y in {6,5,3.5,3,1.5}
\draw (\yy,\y) node[above]{\footnotesize{$\vdots$}};
\foreach  \yy in {1.75,2.75,6.25}
\foreach  \y in {3.5,3}
\draw (\yy,\y) node[above]{\footnotesize{$\vdots$}};
\foreach  \yy in {1.75,2.75,4.25,4.75,6.25}
\foreach  \y in {6.55,5.55,1.05}
\draw (\yy,\y) node[above]{\footnotesize{$\cdots$}};
\foreach  \yy in {4.25,4.75}
\foreach  \y in {6.05,5.05,1.55}
\draw (\yy,\y) node[above]{\footnotesize{$\cdots$}};
\foreach  \yy in {1.75,2.75,6.25}
\foreach  \y in {6,5,1.5}
\draw (\yy,\y) node[above]{\footnotesize{$\ddots$}};

%Indices blocs
\draw[color=blue] (2,1) node[below]{\textbf{$1$}};
\draw[color=blue] (4.5,1) node[below]{\textbf{$\el$}};
\draw[color=blue] (6.5,1) node[below]{\textbf{$\m$}};
\draw[<->,color=blue] (1,0.5) -- (7,0.5);
\draw[color=blue] (4,0.5) node[below]{\textbf{$\m$}};
\draw[color=blue] (7,6) node[right]{\textbf{$1$}};
\draw[color=blue] (7,3.5) node[right]{\textbf{$\kk$}};
\draw[color=blue] (7,1.5) node[right]{\textbf{$\g$}};
\draw[<->,color=blue] (7.5,1) -- (7.5,7);
\draw[color=blue] (7.5,4) node[right]{\textbf{$\g$}};
\draw[color=blue] (4.5,3.5) node[fill=white]{\footnotesize{\textbf{$(\kk,\el)$}}};
\end{tikzpicture}&
%%%%%%%%%%%%%%%%%%%%%%
\hspace{-1cm}\begin{tikzpicture}[scale=0.85]

%Bloc bleu
\foreach  \y in {2,3,4,5}
\draw[ultra thick,color=red] (\y,4.75) -- (\y,6);
\foreach  \y in {2,3,4,5}
\draw[ultra thick,color=red] (\y,1) -- (\y,2.25);
\foreach  \y in {2,3,4,5}
\draw[ultra thick,color=red] (\y,2.75) -- (\y,4.25);
\foreach  \y in {2,3,4,5}
\draw[ultra thick,color=red] (1,\y) -- (2.25,\y);
\foreach  \y in {2,3,4,5}
\draw[ultra thick,color=red] (2.75,\y) -- (4.25,\y);
\foreach  \y in {2,3,4,5}
\draw[ultra thick,color=red] (4.75,\y) -- (6,\y);
%Pointillé
\foreach  \y in {2,3,4,5}
\draw[dotted,ultra thick,color=red] (\y,2.25) -- (\y,2.75);
\foreach  \y in {2,3,4,5}
\draw[dotted,ultra thick,color=red] (\y,4.25) -- (\y,4.75);
\foreach  \y in {2,3,4,5}
\draw[dotted,ultra thick,color=red] (2.25,\y) -- (2.75,\y);
\foreach  \y in {2,3,4,5}
\draw[dotted,ultra thick,color=red] (4.25,\y) -- (4.75,\y);

%Bord
\foreach  \y in {1,6}
\draw[ultra thick,color=black] (\y,4.75) -- (\y,6);
\foreach  \y in {1,6}
\draw[ultra thick,color=black] (\y,1) -- (\y,2.25);
\foreach  \y in {1,6}
\draw[ultra thick,color=black] (\y,2.75) -- (\y,4.25);
\foreach  \y in {1,6}
\draw[ultra thick,color=black] (1,\y) -- (2.25,\y);
\foreach  \y in {1,6}
\draw[ultra thick,color=black] (2.75,\y) -- (4.25,\y);
\foreach  \y in {1,6}
\draw[ultra thick,color=black] (4.75,\y) -- (6,\y);
%Pointillé
\foreach  \y in {1,6}
\draw[dotted,ultra thick,color=black] (\y,2.25) -- (\y,2.75);
\foreach  \y in {1,6}
\draw[dotted,ultra thick,color=black] (\y,4.25) -- (\y,4.75);
\foreach  \y in {1,6}
\draw[dotted,ultra thick,color=black] (2.25,\y) -- (2.75,\y);
\foreach  \y in {1,6}
\draw[dotted,ultra thick,color=black] (4.25,\y) -- (4.75,\y);

\draw (0,0) node{ };

%Alpha
\draw[color=red] (1.5,5.5) node{\textbf{$\al_{11}$}};
\draw[color=red] (3.5,5.5) node{\textbf{$\al_{1\el}$}};
\draw[color=red] (5.5,5.5) node{\textbf{$\al_{1\m}$}};
\draw[color=red] (1.5,3.5) node{\textbf{$\al_{\kk1}$}};
\draw[color=red] (3.5,3.5) node{\textbf{$\alkl$}};
\draw[color=red] (5.5,3.5) node{\textbf{$\al_{\kk\m}$}};
\draw[color=red] (1.5,1.5) node{\textbf{$\al_{\g1}$}};
\draw[color=red] (3.5,1.5) node{\textbf{$\al_{\g\el}$}};
\draw[color=red] (5.5,1.5) node{\textbf{$\al_{\g\m}$}};

%Rho
\draw[color=red] (1.5,1) node[below]{\textbf{$\rhoo_{1}$}};
\draw[color=red] (3.5,1) node[below]{\textbf{$\rhol$}};
\draw[color=red] (5.5,1) node[below]{\textbf{$\rhoo_{\m}$}};

%pi
\draw[color=red] (6,5.5) node[right]{\textbf{$\pii_{1}$}};
\draw[color=red] (6,3.5) node[right]{\textbf{$\pik$}};
\draw[color=red] (6,1.5) node[right]{\textbf{$\pii_{\g}$}};

\end{tikzpicture}\\
\end{tabular}
\end{center}

%% file: ejs-Theorem1Propositions.tex
\section{Proofs of section \ref{sec:profile-likelihood}}
\label{sec:proofPL}
%\subsection{Proof of Proposition~\ref{prop:LocalAsymp}} 

\subsection{Proof of Proposition~\ref{prop:profile-likelihood} (maximum of $\G$ and $\Lambtilde$ in $\theta$)} 
\proofbegin
Define
    $\nu(x, \al)  = x \al - \norm (\al)$. For $\x$ fixed, $\nu(x, \al)$ is maximized at $\al = \normpm(x)$. Manipulations yield
\begin{align*}
  &\Fnd(\bal, \bz, \bw)= \log \prob(\bx; \bz,\bw,\btheta) - \log {\prob(\bx; \bzs,\bws,\bthetas)} \\
%   & = \n\dd \left[ \sum_{\kk} \sum_{\el} \pizk \rhowl \left\{\hxklzw \alkl - \norm(\alkl) \right\}  - \sum_{\kk} \sum_{\el} \pizsk \rhowsl \left\{\hxkl(\bzs, \bws) \alskl - \norm(\alskl) \right\} \right] \\
  & = \n\dd \left[ \sum_{\kk} \sum_{\el} \pizk \rhowl \nu(\hxklzw, \alkl)  - \sum_{\kk} \sum_{\el} \pizsk \rhowsl \nu(\hxkl(\bzs, \bws), \alskl) \right]
\end{align*}
which is maximized at $\alkl = \normpm(\hxklzw)$. Similarly
\begin{align*}
 & \G(\bal, \bz, \bw)  = \Esp_{\bthetas} [ \log \prob(\bx; \bz,\bw,\btheta) - 
\log \prob(\bx; \bzs,\bws,\bthetas) | \bzs, \bws] \\
  & = \n\dd \left[ \sum_{\kk} \sum_{\el} \pizk \rhowl \nu(\barxklzw, \alkl)  - \sum_{\kk} \sum_{\el} \pizsk \rhowsl \nu(\normp(\alskl), \alskl) \right]
\end{align*}
is maximized at $\alkl = \normpm(\barxklzw)$
\proofend

\subsection{Proof of Proposition~\ref{prop:maximum-conditional-likelihood}  (maximum of $\G$ and $\Lambtilde$  in $(\theta,\bz,\bw)$)} 

\proofbegin
We condition on $(\bzs, \bws)$ and prove Equation~\eqref{eq:conditional-likelihood-second-form}:
\begin{align*}
  \G(\btheta, \bz, \bw) & = \Esp_{\bthetas} \left[ \left. \log \frac{\prob(\bx; 
\bz,\bw,\btheta)}{\prob(\bx; \bzs,\bws,\bthetas)} \right| \bzs, \bws \right] \\
  & = \sum_{\ii} \sum_{\jj} \sum_{\kk,\kp}\ \sum_{\el,\lp} \Esp_{\bthetas}\left[ x_{ij} (\al_{\kp\lp} - \alskl) - (\norm(\al_{\kp\lp}) - \norm(\alskl)) \right] \z^\vrai_{\ii\kk} \z_{\ii\kp} \w^\vrai_{\jj\el} \w_{\jj\lp} \\
  & = \n\dd \sum_{\kk, \kp} \sum_{\el, \lp}  \Rgbz_{\kk,\kp} \Rmbw_{\el, \lp} \left[ \normp(\alskl) (\al_{\kp\lp} - \alskl) +  \norm(\alskl) - \norm(\al_{\kp\lp}) \right] \\
  & = - \n\dd \sum_{\kk,\kp} \sum_{\el,\lp} \Rgbz_{\kk,\kp} \Rmbw_{\el, \lp} \KL(\al^\vrai_{\kk\el}, \al_{\kp\lp})
\end{align*}

If $(\bzs,\bws)$ is regular, and for $n, d > 2/c$, all the rows of $\Rgbz$ and $\Rmbw$ have
at least one positive element and we can apply
lemma~\ref{lemme:casdegalite} (which is an adaptation for LBM of Lemma 3.2 of \cite{bickel2013asymptotic} for SBM) to characterize the maximum for
$\G$.

The maximality  of $\Lambtilde(\bzs, \bws)$ results from the fact that $\Lambtilde(\bz, \bw) = \G(\bar{\bal}(\bz, \bw), \bz, \bw)$ where
$\bar{\bal}(\bz, \bw)$ is a particular value of $\bal$, $\Lambtilde$ is
immediately maximum at $(\bz, \bw) \sim (\bzs,\bws)$, and for those, we have $\bar{\bal}(\bz, \bw) \sim \bals$.

The separation and local behavior of $G$ around $(\bzs, \bws)$ is a direct consequence of the proposition \ref{prop:profile-likelihood-derivative}.

\proofend

\subsection{Proof of Proposition~\ref{prop:profile-likelihood-derivative} (Local upper bound for  $\Lambtilde$)}

\proofbegin

We work conditionally on $(\bzs, \bws)$. The principle of the proof relies on the extension of $\Lambtilde$ to a continuous subspace of $M_g([0, 1]) \times M_m([0, 1])$, in which confusion matrices are naturally embedded. The regularity assumption allows us to work on a subspace that is bounded away from the borders of $M_g([0, 1]) \times M_m([0, 1])$. The proof then proceeds by (1) computing the gradient of $\Lambtilde$ at and around its argmax and (2) using those gradients to control the local behavior of $\Lambtilde$ around its argmax. The local behavior allows  in turn to show that $\Lambtilde$ is well-separated.

 Note that $\Lambtilde$ only depends on $\bz$ and $\bw$ through $\Rgbz$ and $\Rmbw$. We can therefore extend it to matrices $(U, V) \in \mathcal{U}_c \times\mathcal{V}_c$ where $\mathcal{U}$ is the subset of matrices $\mathcal{M}_{\g}([0, 1])$ with each row sum higher than $c/2$ and $\mathcal{V}$ is a similar subset of $\mathcal{M}_{\m}([0, 1])$. 
\[
  \Lambtilde(U, V)  = - \n\dd \sum_{\kk,\kp} \sum_{\el,\lp} U_{\kk\kp} V_{\el\lp} \KL\left(\alskl, \baral_{\kp\lp} \right)
  \] where
\[ \baralkl = \baralkl(U, V) = \normpm \left( \frac{\left[U\transpose \Sal V \right]_{\kk\el}}{\left[U\transpose \mathbf{1} V \right]_{\kk\el}} \right) 
\]
and $\mathbf{1}$ is the $\g \times \m$ matrix filled with $1$. Confusion matrices $\Rgbz$ and $\Rmbw$ satisfy $\Rgbz \un = {\widehat \bpi(\bzs)}$ and $\Rmbw \un = \widehat{\brho}(\bws)$, with $\un=(1,\ldots, 1)\transpose$ a vector only containing $1$ values, and are obviously in $\mathcal{U}_c$ and $\mathcal{V}_c$ as soon as $(\bzs, \bws)$ is $c/2$ regular. 

The maps $f_{\kk,\el}: (U, V) \mapsto KL(\alskl, \baralkl(U, V))$ are twice differentiable with second derivatives bounded over $\mathcal{U}_c \times\mathcal{V}_c$ and therefore so is $\Lambtilde(U, V)$. Tedious but straightforward computations show that the derivative of $\Lambtilde$ at $(D_{\pii}, D_{\rhoo}) \coloneqq (\Diag(\widehat{\bpi}(\bzs)), \Diag(\widehat{\brho}(\bws)))$ is: 
\begin{align*}
A_{\kk\kp}(\bws) & \coloneqq -\frac{1}{\n\dd} \frac{\partial 
\Lambtilde}{\partial U_{\kk\kp}}(D_{\pii}, D_{\rhoo}) = \sum_{\el} \widehat{\rhol}(\bws) 
\KL\left(\alskl, \als_{\kp\el} \right) \\
B_{\el\lp}(\bzs) & \coloneqq -\frac{1}{\n\dd} \frac{\partial 
\Lambtilde}{\partial V_{\el\lp}}(D_{\pii}, D_{\rhoo}) = \sum_{\kk} \widehat{\pik}(\bzs) 
\KL\left(\alskl, \als_{\kk\lp} \right)
\end{align*}
$A(\bws)$ and $B(\bzs)$ are the matrix-derivative of $-\Lambtilde/\n\dd$ at $(D_{\pii}, D_{\rhoo})$. Since $(\bzs, \bws)$ is $c/2$-regular and by definition of $\delta(\bals)$, $A(\bws)_{\kk\kp} \geq c\delta(\bals)/2$ (resp. $B(\bws)_{\el\lp} \geq c\delta(\bals)/2$) if $\kk \neq \kp$ (resp. $\el \neq \lp$) and $A(\bws)_{\kk\kk} = 0$ (resp. $B(\bzs)_{\el\el} = 0$) for all $\kk$ (resp. $\el$).   By boundedness of the second derivative, there exists $C > 0$ such that for all $(D_{\pii}, D_{\rhoo})$ and all $(H, G) \in S(D_{\pii}, D_{\rhoo}, C)$, where the definition of the set of local assignments $S$ is extended to the subset of matrices, we have:
\begin{align*}
\frac{-1}{nd} \frac{\partial \Lambtilde}{\partial U_{\kk\kp}}(H, G) \begin{cases} \geq \frac{3c\delta(\bals)}{8} \text { if } \kk \neq \kp \\ \leq \frac{c\delta(\bals)}{8} \text { if } \kk = \kp \end{cases} \text{ and } \quad 
\frac{-1}{nd} \frac{\partial \Lambtilde}{\partial V_{\el\lp}}(H, G) \begin{cases} \geq \frac{3c\delta(\bals)}{8} \text { if } \el \neq \lp \\ \leq \frac{c\delta(\bals)}{8} \text { if } \el = \lp \end{cases}
\end{align*}
Choose $U$ and $V$ in $(\mathcal{U}_{c} \times \mathcal{V}_{c}) \cap S(D_{\pii}, D_{\rhoo}, C)$ satisfying $U\un = \bpi(\bzs)$ and $V\un = \brho(\bws)$. $U - D_{\pii}$ and $V - D_{\rhoo}$ have nonnegative off diagonal coefficients and negative diagonal coefficients. Furthermore, the coefficients of $U, V, D_{\pii}, D_{\rhoo}$ sum up to $1$ and $\Trace(D_{\pii}) = \Trace(D_{\rhoo}) = 1$. By Taylor expansion, there exists a couple $(H, G)$ also in $(\mathcal{U}_c \times \mathcal{V}_c) \cap S(D_{\pii}, D_{\rhoo}, C)$ such that
\begin{multline*}
\frac{-1}{nd} \Lambtilde\left( U, V \right) =  \frac{-1}{nd} \left[\Lambtilde\left( D_{\pii}, D_{\rhoo} \right) + \Trace\left((U - D_{\pii}) \frac{\partial \Lambtilde}{\partial U}(H, G) \right) \right. \\
 \left.  + \Trace\left((V - D_{\rhoo}) \frac{\partial \Lambtilde}{\partial V}(H, G) \right)
\right] \\
\geq \frac{c\delta(\bals)}{8} \left[ 3 \sum_{\kk \neq \kp} (U - D_{\pii})_{\kk\kp} + 3 \sum_{\el \neq \lp} (V-D_{\rhoo})_{\el\lp} + \sum_{\kk} (U - D_{\pii})_{\kk\kk} + \sum_{\el} (V-D_{\rhoo})_{\el\el} \right] \\
= \frac{c\delta(\bals)}{4} [ (1 - \Trace(U)) + (1 - \Trace(V))]
\end{multline*}
To conclude the proof, assume without loss of generality that $(\bz, \bw) \in S(\bzs, \bws, C)$ achieves the $\|.\|_{0,\sim}$ norm (i.e. it is the closest to $(\bzs, \bws)$ in its representative class). Then $(U, V) = (\Rgbz, \Rmbw)$ is in $(\mathcal{U}_c \times \mathcal{V}_c) \cap S(D_{\pii}, D_{\rhoo}, C)$ and satisfy $U\un = \bpi(\bzs)$ (resp. $V\un = \brho(\bws)$). We just need to note $\n(1 - \Trace(\Rgbz)) = \| \bz - \bzs\|_{0,\sim}$ (resp. $\dd(1 - \Trace(\Rmbw)) = \| \bw - \bws\|_{0,\sim}$) to end the proof. \proofend

\section{Proofs of section \ref{sec:big-theorem}}\label{sec:proofCVG}
\subsection{Proof of Proposition~\ref{prop:conditional-likelihood-convergence} (large deviation for $\Fnd$)}

\proofbegin
Conditionally upon $(\bzs, \bws)$,
\begin{align*}
   \Fnd(\btheta, \bz, \bw) - \Lambtilde(\bz,\bw) & \leq \Fnd(\btheta, \bz, \bw) - \G(\btheta, \bz,\bw) \\
  & = \sum_{\ii} \sum_{\jj} (\al_{\zi\wj} - \al^\vrai_{\zi^\vrai \wj^\vrai}) \left(x_{ij} - \normp(\al^\vrai_{\zi^\vrai \wj^\vrai}) \right)  \\
  & = \sum_{\kk \kp} \sum_{\el \lp} \left(\al_{\kp\lp} - \al^\vrai_{\kk \el} \right) W_{\kk \kp \el \lp} \\
  & \leq \sup_{\substack{\Gamma \in \R^{\g^2 \times \m^2} \\ \| \Gamma \|_{\infty} \leq \Diam(\bTheta)}} \sum_{\kk \kp} \sum_{\el \lp} \Gamma_{\kk \kp \el \lp} W_{\kk \kp \el \lp} \coloneqq Z
\end{align*}
uniformly in $\btheta$, where the $W_{\kk \kp \el \lp}$ are independent and defined by:
\begin{equation*}
  W_{\kk \kp \el \lp} = \sum_{\ii} \sum_{\jj} \z^\vrai_{\ii \kk} \w^\vrai_{\jj \el} \z_{\ii, \kp} \w_{\jj \lp}\left(x_{ij} - \normp(\al^\vrai_{\kk\el}) \right)
\end{equation*}
is the sum of $\n\dd \Rgbz_{\kk\kp} \Rmbw_{\el\lp}$ sub-exponential variables with parameters $(\vmax^2, 1/\neighborsize)$ and is therefore itself sub-exponential with parameters \linebreak$(\n\dd \Rgbz_{\kk\kp} \Rmbw_{\el\lp} \vmax^2, 1 / \neighborsize)$. According to Proposition~\ref{prop:concentration-subexponential}, \linebreak $\Esp_{\bthetas}[Z|\bzs, \bws] \leq \g \m \Diam(\bTheta) \sqrt{nd \vmax^2}$ and 
$Z$ is sub-exponential with parameters $(\n\dd \Diam(\bTheta)^2 
(2\sqrt{2})^2\vmax^2, 2\sqrt{2}\Diam(\bTheta) / \neighborsize)$. In particular, 
for all $\vareps_{\n,\dd} < \vmax \neighborsize$
\begin{multline*}
  \Prob_{\bthetas}\left( \left. Z \geq \vmax \g \m \Diam(\bTheta) \sqrt{nd} \left\{ 1 + \frac{\sqrt{8\n\dd} \vareps_{\n,\dd}}{\g \m} \right\} \right| \bzs, \bws \right) \\
  \leq  \Prob_{\bthetas}\left( \left. Z \geq \Esp_{\bthetas}[Z|\bzs, \bws] + \vmax \Diam(\bTheta) \n\dd 2\sqrt{2}\vareps_{\n,\dd} \right| \bzs, \bws \right) \\  \leq \exp \left( - \frac{\n\dd\vareps^2_{\n,\dd}}{2} \right)
\end{multline*}
We can then remove the conditioning and take a union bound to prove Equation~\eqref{eq:conditional-likelihood-convergence}. 
\proofend

\subsection{Proof of Proposition~\ref{prop:large-deviations-profile-likelihood} (contribution of far away assignments)}

\proofbegin
Conditionally on $(\bzs, \bws)$, we know from
proposition~\ref{prop:maximum-conditional-likelihood} that
$\Lambtilde$ is maximal in $(\bzs, \bws)$ and its equivalence
class. Choose $0 < \tnd$ decreasing to $0$ but satisfying $\tnd^2 \gg 
\frac{\n+\dd}{\n\dd}$. This 
is possible as $\n \to 0$ and $\dd \to 0$.

According to equation \ref{eq:conditional-likelihood-upperbound}, for all $(\bz, \bw) \notin 
S(\bzs, \bws, C)$
\[
\Lambtilde(\bz, \bw) \leq - B(C)\; \n\dd.
\]

Now, according to Equation
\ref{eq:conditional-likelihood-separability}, for all $(\bz, \bw) \in S(\bzs, \bws, C) \setminus
S(\bzs, \bws, \tnd)$
\begin{equation*}
\Lambtilde(\bz, \bw) \leq - \frac{c \delta(\bals)}{4}(\n \|\bw - \bws\|_{0, \sim} + \dd \|\bz - \bzs\|_{0, \sim}) \leq - \frac{c \delta(\bals)}{4} \n\dd\tnd
\end{equation*}
since either $\|\bz - \bzs\|_{0, \sim} \geq \n \tnd$ or $\|\bw - \bws\|_{0, \sim} \geq \dd \tnd$. 
% 
% By continuity \CK{?? separation cf proposition \ref{prop:maximum-conditional-likelihood} (iii) . et: est-on sur qu'on peut avoir des $\Lambtilde$ de cet ordre de grandeur?} of
% $\Lambtilde$ and since $\Lambtilde(\bzs, \bws) = 0$, there exist a $\tilde{c} > 0$ such that
% \begin{equation*}
% \forall (\bz, \bw) \notin S(\bzs, \bws, \tilde{c}) \quad \Lambtilde(\bz, \bw) \leq - 2 \n \dd \tnd.
% \end{equation*}
% If additionally $\tilde{c} \leq c_3$, using the definition of $S(\bzs, \bws, \tilde{c})$ , we get 
% \begin{equation}
% \label{eqn:maxlambdatilde}
% \sup_{S^{c}(\bzs, \bws, \tilde{c})} \Lambtilde(\bz, \bw) \leq - \frac{c \tilde{c} \delta(\bals)}{4}\n\dd \ll - \n \dd \tnd
% \end{equation}
% thanks to proposition~\ref{prop:profile-likelihood-derivative}. 
% 
% 
% \CK{Attention: si on n'est pas sur $S$, c'est que $\inf_{\bz' \sim \bz} \|\bz' - \bzs\|_0 > \tilde{c} \n$  OU (et on ET) $\inf_{\bw' \sim \bw} \|\bw' - \bws\|_0 \leq \tilde{c} \dd$. et donc, je ne vois pas comment on arrive au facteur en $nd$. }We can therefore choose $\tilde{c} \leq \min(c_3, c/4)$ without loss of generality. 
%
Hence, for $\n$ and $\dd$ large enough for all $(\bz, \bw) \notin S(\bzs, \bws, \tnd)$
\begin{equation}\label{eqn:maxlambdatilde}
\Lambtilde(\bz, \bw) \leq - \frac{c \delta(\bals)}{4} \n\dd\tnd
\end{equation}
Set $\vareps_{n\dd} = \inf \left( \frac{c\delta(\bals) \tnd 
}{32\sqrt{2}\vmax\Diam(\bTheta)}, \neighborsize \vmax\right)$. 
By 
proposition~\ref{prop:conditional-likelihood-convergence}, and with our choice 
of $\vareps_{n\dd}$, with probability higher than $1 - \Delta_{\n\dd}^1(\vareps_{\n\dd})$,

\begin{align*}
  &\sum_{(\bz, \bw) \notin S(\bzs, \bws, \tnd)} \prob(\bx, \bz, \bw; \btheta)\\ 
  & = \prob(\bx| \bzs, \bws, \bthetas) \sum_{(\bz, \bw) \notin S(\bzs, \bws, \tnd)} \prob(\bz, \bw; \btheta) e^{\Fnd(\btheta, \bz, \bw) - \Lambtilde(\bz, \bw) + \Lambtilde(\bz, \bw)} \\
  & \leq \prob(\bx| \bzs, \bws, \bthetas) \sum_{(\bz, \bw)\in\mcZ\times\mcW} \prob(\bz, \bw; \btheta) e^{\Fnd(\btheta, \bz, \bw) - \Lambtilde(\bz, \bw) - \n \dd \tnd c\delta(\bals) / 4} \\
  & \leq \prob(\bx| \bzs, \bws, \bthetas) \sum_{\bz, \bw} \prob(\bz, \bw; 
\btheta) e^{-\n \dd \tnd c\delta(\bals) / 8} \\
  & = \frac{\prob(\bx, \bzs, \bws; \bthetas)}{\prob(\bzs, \bws; \bthetas)} 
e^{-\n \dd \tnd c\delta(\bals) / 8} \\
  & \leq \prob(\bx, \bzs, \bws; \bthetas) \exp\left( -\n\dd \tnd
\frac{c\delta(\bals)}{8} + (\n+\dd)\log \frac{1}{c} \right) \\
  & = \prob(\bx, \bzs, \bws; \bthetas) \smallO(1)
\end{align*}
where the second line comes from inequality (\ref{eqn:maxlambdatilde}), the third from the global control studied in Proposition~\ref{prop:conditional-likelihood-convergence} and the definition of $\vareps_{\n\dd}$, the fourth from the definition of $\prob(\bx, \bzs, \bws; \bthetas)$, the fifth from the bounds on $\bpis$ and $\brhos$ and the last from $\tnd \gg (\n+\dd)/\n\dd$. 

In addition, we have $\vareps_{n\dd}^2 \gg \frac{\n+\dd}{\n\dd}$ so that $\Delta_{\n\dd}^1(\vareps_{\n\dd})$ vanishes and:
\begin{align*}
  \sum_{(\bz, \bw) \notin S(\bzs, \bws, \tnd)} \prob(\bx, \bz, \bw; \btheta) & = \prob(\bx, \bzs, \bws, \bthetas) \smallO_P(1)
\end{align*}
where the $\smallO_P$ is uniform in probability over all $\btheta \in\bTheta$.
\proofend
  
\subsection{Proof of Proposition~\ref{prop:profile-likelihood-convergence-local} (local convergence $\Fnd$)} 

\proofbegin
We work conditionally on $(\bzs, \bws) \in \mcZ_1 \times \mcW_1$. We assume with no loss of generality that $(\bz,\bw)$ is the representative of its class closest to $(\bzs,\bws)$. Choose
$\vareps \leq \neighborsize \vmin^2$ small.  Manipulation of $\Lamb$ and $\Lambtilde$ yield
\begin{align*}
  \frac{\Fnd(\btheta, \bz, \bw) - \Lambtilde(\bz, \bw)}{\n\dd} \leq &
\frac{\Lamb(\bz, \bw) - \Lambtilde(\bz, \bw)}{\n\dd} \\
   = & \sum_{\kk} \sum_{\el} \pizk\rhowl \left[ 
f(\hx_{\kk\el}) - f(\barx_{\kk\el}) \right] \\
   & - \sum_{\kk} \sum_{\el} \pizsk \rhowsl \alskl 
 (\hxkl^\star - \barxkl^\star) 
\end{align*}
where $f(x) =  x\normpm(x) - \norm\circ\normpm(x)$, $\hxkl^\star = 
\hxkl(\bzs, \bws)$ and $\barxkl^\star = \normp(\alskl)$. The 
function $f$ is twice differentiable on $\mathring{\mathcal{A}}$ with $f'(x) = 
-\normpm(x)$ and $f''(x) = - 1 / \norm'' \circ \normpm(x)$. $f'$ (resp. $f''$)
are bounded over $I = \normp([-M_{\al}, M_{\al}])$ by $M_{\al}$ 
(resp. $1/\vmin^2$).

Assignments $(\bz, \bw)$ %at $\|.\|_{0,\sim}$-distance less than $c/4$ of $(\bzs, \bws)$ 
belonging to $S(\bzs, \bws,c/4)$ are also $c/4$-regular. 
According to Proposition~\ref{proposition:maxzw}, $\hxkl$ and $\barxkl$ are 
at distance at most $\vareps$ with probability higher than $1 - 2 \exp \left( - 
\frac{\n\dd c^2 \vareps^2}{32(\vmax^2 + \neighborsize^{-1} 
\vareps)}\right)$, so that:
\begin{equation*}
  f(\hxkl) - f(\barxkl) = f'(\barxkl) \left( \hxkl - \barxkl \right) + 
\Omega\left( (\hxkl - \barxkl)^2 \right)
\end{equation*}
By Proposition~\ref{proposition:maxzw}, $(\hxkl - \barxkl)^2 = 
\bigO_P(1/\n\dd)$ where the $\bigO_P$ is uniform in $\bz, \bw$ and does not 
depend on $\bzs, \bws$. Similarly, 
\begin{equation*}
  f'(\barxkl) = f'(\barxkl^\star) + \Omega(\barxkl - \barxkl^\star) = \alskl 
+ \Omega(\barxkl - \barxkl^\star)
\end{equation*}
$\barxkl$ is a convex combination of the $\Salkl = \normp(\alskl)$ therefore, 
\begin{align*}
  | \barxkl - \barxkl^\star | &= \left| \frac{\left[\Rgbz\transpose \Sal \Rmbw 
\right]_{\kk\el}}{\pizk\rhowl} - \barxkl^\star \right| \\
  & \leq \left( 1 - \frac{\Rgbz_{\kk\kk}\Rmbw_{\el\el}}{\pizk\rhowl}\right) 
(\Salmax - \Salmin)
\end{align*}
Note that: 
\begin{align*}
  \sum_{\kk,\el} \pizk \rhowl \left( 1 - 
\frac{\Rgbz_{\kk\kk}\Rmbw_{\el\el}}{\pizk\rhowl}\right)
 &= 1 - \Trace(\Rgbz)\Trace(\Rmbw) \\ 
 & \leq \frac{\|\bz - \bzs\|_{0, \sim}}{\n} + \frac{\|\bw - \bws\|_{0, \sim}}{\dd}
\end{align*}
and $\hxkl - \barxkl = \smallO_P(1)$. Therefore 
\begin{equation*}
\sum_{\kk,\el} \pizk \rhowl \Omega(\barxkl - \barxkl^\star) \times (\hxkl 
- \barxkl) = \smallO_P\left( \frac{\|\bz - \bzs\|_{0, \sim}}{\n} + \frac{\|\bw 
- \bws\|_{0, \sim}}{\dd}\right)
\end{equation*}
The remaining term writes
\begin{align*}
\sum_{\kk,\el} \alskl \left[ \pizk \rhowl (\hxkl - \barxkl) - \pizsk \rhowsl 
(\hxkl^\star - \barxkl^\star)  \right]
\end{align*}
According to Proposition~\ref{proposition:maxdiffzw}, this term is $\smallO_P(\left( \frac{\|\bz - \bzs\|_{0, \sim}}{\n}  + \frac{\|\bw - 
\bws\|_{0, \sim}}{\dd}\right)$ uniformly in $(\bz, \bw)$ and $(\bzs, \bws) \in 
\Omega_1$ as soon $\n e^{-a\dd} \to 0$ and $\dd e^{-a\n} \to 0$ for all $a > 0$, which is true under $(H_4)$. It follows that:
\begin{equation*}
\sup_{\substack{(\bz, \bw) \nsim (\bzs, \bws) \\ (\bz, \bw) \in S(\bzs, \bws, c/4)}} \frac{\Lamb(\bz, \bw) - 
\Lambtilde(\bzs, \bws)}{\n\dd} = \smallO_P \left( \frac{\|\bz - \bzs\|_{0, \sim}}{\n} + \frac{\|\bw - \bws\|_{0, \sim}}{\dd} \right)
\end{equation*}
\proofend

\subsection{Proof of Proposition~\ref{prop:small-deviations-profile-likelihood} (contribution of local assignments)} 

\proofbegin
By Proposition~\ref{cor:prob-regular-configurations-star}, it is enough to prove 
that the sum is small compared to $\prob(\bzs, \bws, \bx; \bthetas)$ on $\Om_1$. 
We work conditionally on $(\bzs, \bws) \in \mcZ_1 \times \mcW_1$. Choose $(\bz, 
\bw)$ in $S(\bzs, \bws, \tilde{c})$. This set is non empty as soon as $\min(\tilde{c}\n,\tilde{c}\dd) >1$.
\begin{align*}
  \log \left( \frac{\prob(\bz, \bw, \bx; \btheta)}{\prob(\bzs, \bws, \bx; 
\bthetas)} \right) & = \log \left( \frac{\prob(\bz, \bw; \btheta)}{\prob(\bzs, 
\bws; \bthetas)}\right) + \Fnd(\btheta, \bz, \bw) \\
\end{align*}
We can assume without loss of generality that $(\bz, \bw)$ 
is the representative closest to $(\bzs, \bws)$ and denote $r_1 = \|\bz - 
\bzs\|_0$ and $r_2 = \|\bw - \bws\|_0$. Then:
\begin{align*}
  \Fnd(\btheta, \bz, \bw) & \leq \Lamb(\bz, \bw) - \Lambtilde(\bz, \bw) + 
\Lambtilde(\bz, \bw) \\
  & \leq \Lamb(\bz, \bw) - \Lambtilde(\bz, \bw) - 
\frac{c\delta(\bals)}{4}\left(\dd r_1 + \n r_2\right) \\
  & \leq -\frac{c\delta(\bals)}{4}\left(\dd r_1 + \n r_2\right)(1 + 
\smallO_P(1)) \\
\end{align*}
where the first line comes from the definition of $\Lamb$, the second line from 
Proposition~\ref{prop:profile-likelihood-derivative} and the fact that $\tilde{c} < C$ and the third from Proposition~\ref{prop:profile-likelihood-convergence-local} and the fact that $\tilde{c} < c/4$. Thanks to corollary~\ref{cor:marginalprobabilties}, we also know that:
\begin{equation*}
 \log \left( \frac{\prob(\bz, \bw; \btheta)}{\prob(\bzs, \bws; \bthetas)}\right) 
\leq \bigO_P(1) \exp\left\{M_{c/4}(r_1 + r_2) \right\}
\end{equation*}
There are at most ${\n \choose r_1}{\n \choose r_2}\g^{r_1}\m^{r_2}$ assignments 
$(\bz,\bw)$ at distance $r_1$ and $r_2$ of $(\bzs, \bws)$ and each of them has 
at most $\g^{\g} \m^{\m}$ equivalent configurations. Therefore, 
\begin{align*}
 & \frac{\sum_{\substack{(\bz, \bw) \in S(\bzs, \bws, \tilde{c}) \\ (\bz, \bw) 
\nsim (\bzs, \bws)}} \prob(\bz, \bw, \bx; \btheta)}{\prob(\bzs, \bws, \bx; 
\bthetas)} \\
  & \leq \bigO_P(1) \sum_{\substack{r_1 + r_2 \geq 1}} {\n \choose r_1}{\n 
\choose r_2}\g^{\g + r_1}\m^{\m + r_2} \exp\left( (r_1 + r_2)M_{c/4} - 
\frac{c\delta(\bals)}{4}\left(\dd r_1 + \n r_2\right)(1 + \smallO_P(1)) \right) 
\\
  & = \bigO_P(1) \left(1+ e^{\log \g + M_{c/4} - \dd \frac{c 
\delta(\bals)(1 + \smallO_P(1))}{4} }\right)^{\n}
  \left(1+ e^{\log \m + M_{c/4} - \n \frac{c \delta(\bals)(1 + 
\smallO_P(1))}{4} }\right)^{\dd} -1 \\
  & \leq \bigO_P(1) a_{\n\dd} \exp(a_{\n\dd})
\end{align*}
where $a_{\n\dd} = \n e^{\log \g + M_{c/4} - \dd \frac{c 
\delta(\bals)(1 + \smallO_P(1))}{4} } + \dd e^{\log \m + M_{c/4} - \n 
\frac{c \delta(\bals)(1 + \smallO_P(1))}{4} } = \smallO_P(1)$ as soon as $\n 
\gg \log \dd$ and $\dd \gg \log \n$.
\proofend

\subsection{Proof of Proposition~\ref{prop:equivalent-configurations-profile-likelihood} (contribution of equivalent assignments)} 

\proofbegin
Choose $(s, t)$ permutations of $\{1, \dots, \g\}$ and $\{1, \dots, \m\}$ and assume that $\bz = \bz^{\vrai,s}$ and $\bw = \bw^{\vrai, t}$. Then $\prob(\bx, \bz, \bw; \btheta) = \prob(\bx, \bz^{\vrai, s}, \bw^{\vrai, t}; \btheta) = \prob(\bx, \bzs, \bws; \btheta^{s,t})$.  If furthermore $(s, t) \in \Symmetric(\btheta)$, $\btheta^{s, t} = \btheta$ and immediately $\prob(\bx, \bz, \bw; \btheta) = \prob(\bx, \bzs, \bws; \btheta)$. We can therefore partition the sum as 

\begin{align*}
  \sum_{(\bz, \bw) \sim (\bzs, \bws)} \prob(\bx, \bz, \bw; \btheta) & = \sum_{s, t} \prob(\bx, \bz^{\vrai, s}, \bw^{\vrai, t}; \btheta) \\ 
  & = \sum_{s, t} \prob(\bx, \bzs, \bws; \btheta^{s,t}) \\ 
  & = \sum_{\btheta' \sim \btheta} \# \Symmetric(\btheta') \prob(\bx, \bzs, \bws; \btheta') \\ 
  & = \# \Symmetric(\btheta) \sum_{\btheta' \sim \btheta} \prob(\bx, \bzs, \bws; \btheta') \\ 
\end{align*}

The complete likelihood $\prob(\bx, \bzs, \bws; \btheta)$ is a unimodal function of  $\btheta$ with  mode located in $\bthetaMC$. By consistency of $\bthetaMC$, either $\prob(\bx, \bzs, \bws; \btheta) = \smallO_P(\prob(\bx, \bzs, \bws; \bthetas))$ or $\prob(\bx, \bzs, \bws; \btheta) = \bigO_P(\prob(\bx, \bzs, \bws; \bthetas))$ when $\btheta$ is in a close neighborhood of $\bthetas$.                                                                                                                                                %, $\btheta  \to \bthetas$. 
In the latter case, any $\btheta' \sim \btheta$ other than $\btheta$ is bounded away from $\bthetas$ and thus $\prob(\bx, \bzs, \bws; \btheta') = \smallO_P(\prob(\bx, \bzs, \bws; \bthetas))$. In summary, 
\[
 \sum_{\btheta' \sim \btheta} \frac{\prob(\bx, \bzs, \bws; \btheta')}{\prob(\bx, \bzs, \bws; \bthetas)} = \max_{\btheta' \sim \btheta} \frac{\prob(\bx, \bzs, \bws; \btheta')}{\prob(\bx, \bzs, \bws; \bthetas)} (1 + \smallO_P(1))
\] \proofend

\section{Concentration for sub-exponential variables}
\label{sec:subexp}
%%%%%%%%%%%%%%%%%%%%%%%%%

\begin{comment}
\subsection{Proof of Lemma~\ref{cor:prob-regular-configurations-star}}
\CK{A laisser? ou à mettre dans le texte?}
\proofbegin
Each $\zsumk$ is a sum of $\n$ i.i.d Bernoulli r.v. with parameter $\pik \geq \pii_{\min} \geq c$. A simple Hoeffding bound shows that
\begin{equation*}
  \Prob_{\bthetas}\left( \zsumk \leq \n \frac{c}{2} \right)
  \leq
  \Prob_{\bthetas}\left( \zsumk \leq \n \frac{\pik}{2} \right)
  \leq
  \exp\left( - 2\n\left(\frac{\pik}{2}\right)^2 \right)
  \leq
  \exp\left( - \frac{\n c^2}{2} \right)
\end{equation*}
Similarly,
\begin{equation*}
  \Prob_{\bthetas}\left( \wsuml \leq \dd \frac{c}{2} \right) \leq \exp\left( - \frac{\dd c^2}{2} \right)
\end{equation*}
The proposition follows from a union bound over $\g$ values of $\kk$ and $\m$ values of $\el$.
\proofend
\end{comment}
%\section{Sub-exponential variables}
%\label{sec:subexp}

%% Subexponential variables
Concentration inequalities for sub-exponential variables play a key role: in particular  Proposition~\ref{prop:concentration-subexponential} for global convergence and Propositions~\ref{proposition:maxzw} and \ref{proposition:maxdiffzw} for local convergence. We present here some properties of sub-exponential variables (\cite{wainwright2019high}), then derives the  needed concentration inequalities.

Recall first that a random variable $X$ is  sub-exponential with parameters $(\tau^2, b)$ if for all $\lambda$ such that $|\lambda|\leq 1/b$,
\[
\Esp[e^{\lambda( X-\Esp(X))}] \leq \exp\left( \frac{\lambda^2 \tau^2}{2} \right).
\]
In particular, all distributions coming from a natural exponential family are sub-exponential. Sub-exponential variables satisfy a large deviation Bernstein-type inequality:
\begin{equation*}
  \label{eq:concentration-subexponential-def}
  \Prob( X - \Esp[X] \geq t) \leq 
  \begin{cases} 
    \exp\left( - \frac{t^2}{2 \tau^2}\right)  & \text{if} \quad 0 \leq t \leq \frac{\tau^2}{b} \\ 
    \exp\left( - \frac{t}{2b}\right)  & \text{if} \quad t \geq \frac{\tau^2}{b}
  \end{cases}
\end{equation*}
So that
\[
\Prob( X - \Esp[X] \geq t) \leq \exp\left( - \frac{t^2}{2(\tau^2+bt)}\right)
\]
\subsection{Properties}
The sub-exponential property is preserved by summation and multiplication. 
\begin{itemize}
\item If $X$ is sub-exponential with parameters $(\tau^2, b)$ and $\alpha \in \mathbb{R}$, then so is $\alpha X$ with parameters $(\alpha^2\tau^2, \alpha b)$
\item If the $X_i$, $i = 1,\dots,n$ are sub-exponential with parameters $(\tau_i^2, b_i)$ and independent, then so is $X = X_1 + \dots + X_n$ with parameters $(\sum_i \tau_i^2,\max_i b_i)$
\end{itemize}
Moreover, Lemma \ref{lemme:absolutesubexponential} defines the sub-exponential property of the absolute value of a sub-exponential variable.
\begin{lemme}\label{lemme:absolutesubexponential}\textbf{}
If $X$ is a zero mean random variable, sub-exponential with parameters $(\sigma^2, b)$, then $|X|$ is sub-exponential with parameters $(8\sigma^2, 2\sqrt{2}b)$.
\end{lemme}

\proofbegin
Denote $\mu = \Esp|X|$ and consider $Y = |X| - \mu$. Choose $\lambda$ such that $|\lambda| < (2\sqrt{2}b)^{-1}$. We need to bound $\Esp[e^{\lambda Y}]$. Note first that $\Esp[e^{\lambda Y}] \leq \Esp[e^{\lambda X}] + \Esp[e^{-\lambda X}] < +\infty$ is properly defined by sub-exponential property of $X$ and we have
\[ 
\Esp[e^{\lambda Y}] \leq 1 + \sum_{k=2} \frac{|\lambda|^k\Esp[|Y|^k]}{k!} 
\]
where we used the fact that $\Esp[Y] = 0$. We know bound odd moments of $|\lambda Y|$. 
\[
\Esp[|\lambda Y|^{2k+1}] \leq (\Esp[|\lambda Y|^{2k}]\Esp[|\lambda Y|^{2k+2}])^{1/2} \leq \frac{1}{2} (\lambda^{2k} \Esp[Y^{2k}] + \lambda^{2k+2} \Esp[Y^{2k+2}])
\]
where we used first Cauchy-Schwarz and then the arithmetic-geometric mean inequality. The Taylor series expansion can thus be reduced to 
\begin{align*}
\Esp[e^{\lambda Y}] & \leq 1 + \left(\frac{1}{2} + \frac{1}{2.3!}\right) \Esp[Y^2] \lambda^2+ \sum_{k=2}^{+\infty} \left( \frac{1}{(2k)!} + \frac{1}{2}\left[ \frac{1}{(2k-1)!} + \frac{1}{(2k+1)!} \right] \right)\lambda^{2k} \Esp[Y^{2k}] \\ 
& \leq \sum_{k=0}^{+\infty} 2^k \frac{\lambda^{2k}\Esp[Y^{2k}]}{(2k)!} \\
&\leq \sum_{k=0}^{+\infty} 2^{3k} \frac{\lambda^{2k}\Esp[X^{2k}]}{(2k)!} =\cosh\left(2\sqrt{2}\lambda X\right) 
=\Esp\left[ \frac{e^{2\sqrt{2}\lambda X} + e^{-2\sqrt{2}\lambda X}}{2}\right]\\
& \leq e^{\frac{8\lambda^2\sigma^2}{2}}\\
\end{align*}
where we used the well-known inequality $\Esp[|X - \Esp[X]|^k] \leq 2^k \Esp[|X|^k]$ to substitute $2^{2k}\Esp[X^{2k}]$ to $\Esp[Y^{2k}]$.\\

\proofend

\subsection{Concentration inequalities}

\begin{proposition}[Maximum in $(\bz, \bw)$] \label{proposition:maxzw}
  Let $(\bz, \bw)$ be a configuration and $\widehat{x}_{\kk, \el}(\bz, \bw)$ resp. $\barxkl(\bz, \bw)$ be  defined in Equations~\eqref{eq:mle-complete-likelihood} and~\eqref{eq:profile-likelihood-notations}. Under the assumptions of the section~\ref{sec:assumptions}, for all $\vareps>0$
  \begin{equation}
    \label{eq:ineg-de-bernst-1}
    \Prob \left(\max_{\bz, \bw} \max_{k,l} \pizk \rhowl |\widehat{x}_{\kk, \el} - \barxkl| > \vareps \right)  \leq 2 \g^{\n+1} \m^{\dd + 1} \exp \left( - \frac{nd\vareps^2}{2(\bar{\sigma}^2 + \neighborsize^{-1}\vareps)}\right). \\
  \end{equation}
  Additionally, the suprema over all $c/2$-regular  assignments satisfies:
  \begin{equation}
    \label{eq:ineg-de-bernst-2}
    \Prob \left(\max_{\bz \in \mcZ_1, \bw \in \mcW_1} \max_{k,l} |\widehat{x}_{\kk, \el} - \barxkl| > \vareps \right) \leq  2 \g^{\n+1} \m^{\dd + 1} \exp \left( - \frac{ndc^2\vareps^2}{8(\bar{\sigma}^2 + \neighborsize^{-1}\vareps)}\right). \\
  \end{equation}
\end{proposition}
Note that equations \ref{eq:ineg-de-bernst-1} and \ref{eq:ineg-de-bernst-2} remain valid when replacing $c/2$ by any $\tilde c < c/2
$.
\proofbegin

The random variables $X_{ij}$ are subexponential with parameters $(\vmax^2, 1/\neighborsize)$. Conditionally to $(\bzs, \bws)$, $\zsumk \wsuml(\widehat{x}_{\kk, \el} - \barxkl)$ is a sum of $\zsumk \wsuml$ centered subexponential random variables. By Bernstein's inequality \cite{massart2007concentration}, we therefore have for all $t > 0$
\begin{align*}
\Prob( \zsumk \wsuml | \widehat{x}_{\kk, \el} - \barxkl| \geq t ) & \leq 2 \exp\left( - \frac{t^2}{2(\zsumk \wsuml\bar{\sigma}^2 + \neighborsize^{-1} t)}\right) \\
\end{align*}
In particular, if $t = \n\dd x$,
\begin{align*}
\Prob\left( \pizk \rhowl | \widehat{x}_{\kk, \el} - \barxkl| \geq x \right) & \leq 2 \exp\left( - \frac{ndx^2}{2(\pizk\rhowl\bar{\sigma}^2 + \neighborsize^{-1} x)}\right)\\ &\leq 2 \exp\left( - \frac{ndx^2}{2(\bar{\sigma}^2 + \neighborsize^{-1} x )}\right) \\
\end{align*}
uniformly over $(\bz, \bw)$. Equation~\eqref{eq:ineg-de-bernst-1} then results from a union bound.
Similarly,
\begin{align*}
  \Prob\left( | \widehat{x}_{\kk, \el} - \barxkl| \geq x \right) & = \Prob\left( \pizk \rhowl | \widehat{x}_{\kk, \el} - \barxkl| \geq \pizk \rhowl x \right) \\
  & \leq 2 \exp\left( - \frac{ndx^2\pizk^2\rhowl^2}{2(\pizk\rhowl\bar{\sigma}^2 + \neighborsize^{-1} x \pizk\rhowl)}\right) \\
  & \leq 2 \exp\left( - \frac{ndc^2x^2}{8(\bar{\sigma}^2 + \neighborsize^{-1} x )}\right)
\end{align*}
Where the last inequality comes from the fact that $c/2$-regular assignments satisfy $\pizk \rhowl \geq c^2/4$. Equation~\eqref{eq:ineg-de-bernst-2} then results from a union bound over $\mcZ_1 \times \mcW_1 \subset \mcZ \times \mcW$.
\proofend

\begin{proposition}[Maximum in non equivalent $(\bz, \bw)$] \label{proposition:maxdiffzw}
  Let $(\bar{\bz}, \bar{\bw})$ be any configuration and $(\bz, \bw)$ the 
$\sim$-equivalent configuration that achieves $\|\bz - \bzs\|_0 = \|\bar{\bz} - 
\bzs\|_{0,\sim}$ and $\|\bw - \bws\|_0 = \|\bar{\bw} - 
\bws\|_{0,\sim}$, let $\hxkl = \widehat{x}_{\kk, \el}(\bz, \bw)$ (resp. 
$\barxkl(\bz, \bw)$) and $\hxkl^\star = \widehat{x}_{\kk, \el}(\bzs, 
\bws)$ (resp. $\barxkl^\star = \barxkl(\bzs, \bws)$ = $\normp(\alskl)$) be as 
defined in Equations~\eqref{eq:mle-complete-likelihood} 
and~\eqref{eq:profile-likelihood-notations}. Under the assumptions $(H_1)$ to $(H_3)$, for all $\vareps \leq \neighborsize 
\vmax^2$, 
  \begin{multline*}
    \label{eq:ineg-de-bernst-1}
    \Prob \left( \max_{(\bar{\bz}, \bar{\bw}) \nsim (\bzs, \bws)} \max_{k,l} 
\frac{\n\dd\left[ \pizk \rhowl (\widehat{x}_{\kk, \el} - \barxkl) - \pizsk \rhowsl 
(\hxkl^\star - \barxkl^\star ) \right]}{\n \|\bw - \bws\|_{0} + \dd\|\bz 
- 
\bzs\|_{0}} > \vareps \right) \\
\leq \g^{\g+1} \m^{\m+1} a_{\n\dd}e^{a_{\n\dd}}
  \end{multline*}
where $a_{\n\dd} = \n\g e^{-\frac{\dd\vareps^2}{\vmax^2}} + \m\dd e^{-\frac{\n\vareps^2}{\vmax^2}}$.
\end{proposition}

\proofbegin Denote $r_1 = \|\bz - \bzs\|_{0}$ and $r_2 
= \|\bw - \bws\|_{0}$. The numerator within the $\max$ in the fraction 
can be expanded to 
\begin{equation*}
Z_{\kk\el}(\bz, \bw) = \sum_{i, j} (\zik\wjl - \zik^\star \wjl^\star) (X_{ij} 
- \als_{\zik^\star \wjl^\star})
\end{equation*}
and is thus a sum of at most $N = \n r_2 + \dd r_1$ 
non-null centered sub-exponential random variables with parameters $(\vmax^2, 
1/\neighborsize)$. It is therefore centered sub-exponential with 
parameters $(N \vmax^2, 1/\neighborsize)$. By Bernstein inequality, for all 
$\vareps \leq \neighborsize \vmax^2$ we have
\begin{equation*}
\Prob (Z \geq \vareps (\n r_2 + \dd r_1) ) \leq \exp\left(- 	
\frac{(\n r_2 + \dd r_1)\vareps^2}{2\vmax^2} \right).
\end{equation*} 
There are at most ${\n \choose r_ 1} \g^{r_1} \g^\g$ $\bz$ at $\|.\|_{0,\sim}$ distance 
$r_1$ of $\bzs$ and ${\dd \choose r_2} \m^{r_2} \m^\m$ $\bz$ at $\|.\|_{0,\sim}$ 
distance $r_2$ of $\bws$. An union bound shows that:
\begin{multline*}
\Prob \left( \max_{(\bar{\bz}, \bar{\bw}) \nsim (\bzs, \bws)} \max_{k,l} 
\frac{Z_{\kk\el}(\bz, \bw)}{\n \|\bw - \bws\|_{0} + \dd\|\bz - \bzs\|_{0}} 
 \geq \vareps \right) \\ 
 \leq \sum_{r_1 + r_2 \geq 1} \sum_{\substack{\|\bar{\bz} - \bzs\|_{0, 
\sim}=r_1\\ \|\bar{\bw} - \bws\|_{0, \sim}=r_2}} \g\m \Prob (Z_{\kk\el}(\bz, \bw) \geq 
\vareps (\n r_2 + \dd r_1) ) \\
 \leq \sum_{r_1 + r_2 \geq 1} {\n \choose r_1} {\dd \choose r_2} \g^{\g+1} m^{\m+1} \g^{r_1} \m^{r_2} \exp\left(-(\n r_2 + \dd 
r_1)\vareps^2 / 2\vmax^2\right) \\
 = \g^{\g+1} m^{\m+1} \sum_{r_1 + r_2 \geq 1} \left( \g e^{-\frac{\dd\vareps^2}{\vmax^2}} \right)^{r_1} \left( \m e^{-\frac{\n\vareps^2}{\vmax^2}} \right)^{r_2}
 \leq \g^{\g+1} \m^{\m+1} a_{\n\dd}e^{a_{\n\dd}}
\end{multline*}
where $a_{\n\dd} = \n\g e^{-\frac{\dd\vareps^2}{\vmax^2}} + \m\dd e^{-\frac{\n\vareps^2}{\vmax^2}}$. 
\proofend

% \begin{proposition} \label{cor:maxdiffzw}
% With the notations of \ref{proposition:maxdiffzw} and under assumption $(H_4)$:
% \[
%  \max_{(\bar{\bz}, \bar{\bw}) \nsim (\bzs, \bws)} \max_{k,l} 
% \frac{\n\dd\left[ \pizk \rhowl (\widehat{x}_{\kk, \el} - \barxkl) - \pizsk \rhowsl 
% (\hxkl^\star - \barxkl^\star ) \right]}{\n \|\bw - \bws\|_{0} + \dd\|\bz 
% - 
% \bzs\|_{0}} = \smallO_P(1)
% \]
% \end{proposition}
%  
% \proofbegin
% This is a direct consequence from Proposition~\ref{proposition:maxdiffzw} and the fact that $a_{\n\dd} \to 0$ under $(H_4)$.
% \proofend

\begin{proposition}[concentration for sub-exponential]
\label{prop:concentration-subexponential}
Let $\X_{1}, \dots, \X_{\n}$ be independent zero mean random variables, sub-exponential with parameters $(\sigma_i^2, b_i)$. Denote $V_0^2  = \sum_{\ii} \sigma_i^2$ and $b = \max_{i} b_i$. Then the random variable $Z$ defined by:
\begin{equation*}
  Z = \sup_{\substack{\Gamma \in \R^{\n} \\ \|\Gamma \|_{\infty} \leq M}} \sum_{\ii} \Gamma_i X_{\ii}
\end{equation*}
is also sub-exponential with parameters $(8 M^2 V_0^2, 2\sqrt{2} Mb)$. 
% and the following Bernstein inequality holds:
% \begin{equation*}
%   \Prob( Z - \Esp[Z] \geq t) \leq 
%   \begin{cases} 
%     \exp\left( - \frac{t^2}{2 (M^2 V_0^2+n \log(2))}\right)  & \text{if} \quad 0 \leq t \leq \frac{M^2 V_0^2 + \n \log(2)}{Mb} \\ 
%     \exp\left( - \frac{t}{2Mb}\right)  & \text{if} \quad t \geq \frac{M^2 V_0^2 + \n \log(2)}{Mb}
%   \end{cases}
% \end{equation*}
Moreover $\Esp[Z] \leq M V_0 \sqrt{\n}$ so that for all $t > 0$, 
\begin{equation*}
  \label{eq:concentration-subexponential}
\Prob( Z - M V_0 \sqrt{n} \geq t) \leq \exp \left( - \frac{t^2}{2 (8M^2 V_0^2 + 2\sqrt{2} M b t)}\right)
\end{equation*}

\end{proposition}

\proofbegin
Note first that $Z$ can be simplified to $Z =  M \sum_{\ii} |X_{\ii}|$. We just need to bound bound $\Esp[Z]$. The rest of the proposition results from the fact that the $|X_{\ii}|$ are subexponential $(8\sigma_{\ii}^2, 2\sqrt{2}b_{\ii})$ by Lemma~\ref{lemme:absolutesubexponential} and standard properties of sums of independent rescaled subexponential variables. 
\begin{align*}
  \Esp[Z] & = \Esp \left[ \sup_{\substack{\Gamma \in \R^{\n} \\ \|\Gamma \|_{\infty} \leq M}}  \sum_{\ii} \Gamma_{\ii} X_{\ii} \right]  = \Esp \left[ \sum_{\ii} M |X_{\ii}| \right] \leq M \sum_{\ii} \sqrt{\Esp [ X_{\ii}^2 ]} \\ & 
  = M \sum_{\ii} \sigma_i \leq M \left( \sum_{\ii} 1 \right)^{1/2} \left( \sum_{\ii} \sigma_{\ii}^2 \right)^{1/2} = M V_0\sqrt{\n}
\end{align*}
using Cauchy-Schwarz. 
\proofend

\section{Technical lemmas}
\label{sec:technical-lemma}
Lemma \ref{lemme:casdegalite} is the working horse for proving Proposition~\ref{prop:maximum-conditional-likelihood}.  Corollary \ref{cor:marginalprobabilties} is needed for Theorem \ref{prop:small-deviations-profile-likelihood} and Lemma \ref{lemme:marginalprobabilties} is an intermediate result for Corollary \ref{cor:marginalprobabilties}.

\begin{lemme}\label{lemme:casdegalite}\textbf{}

Let $\etaa$ and $\bar{\etaa}$ be two matrices from $M_{\g \times \m}(\Theta)$ and $f: \Theta \times \Theta \to \mathbb{R_+}$ a positive function, $A$ a (squared) confusion matrix of size $\g$ and $B$ a (squared) confusion matrix of size $\m$. We denote $d_{\kk\el\kk'\el'} = f(\etaa_{\kk\el}, \bar{\etaa}_{\kk'\el'})$. Assume that
\begin{itemize}
\item all the rows of $\etaa$ are distinct;
\item all the columns $\etaa$ are distinct;
\item $f(x,y) = 0 \Leftrightarrow x = y$;
\item each  row of $A$ has a non zero element;
\item each  row of $B$ has a non zero element;
\end{itemize}
and denote
\begin{equation*}
\label{eq:casdegalite}
\Sigma = \sum_{\kk \kk'} \sum_{\el \el'} A_{\kk \kk'} B_{\el \el'} d_{\kk\el\kk'\el'}
\end{equation*}
Then,
\begin{equation*}
  \Sigma = 0 \Leftrightarrow \begin{cases}
    A, B \text{ are permutation matrices } s,t & \\
    \bar{\etaa} = \etaa^{s, t} \text{ i.e. } \forall (\kk, \el), \bar{\etaa}_{\kk\el} = \etaa_{s(\kk) t(\el)} &
  \end{cases}
\end{equation*}
\end{lemme}
\proofbegin
If $A$ and $B$ are the permutation matrices corresponding to the permutations $s$ et $t$: $A_{ij} = 0$ if $i \neq s(j)$ and $B_{ij} = 0$ if $i \neq t(j)$. As each row of $A$ contains a non zero element  and as $A_{s(\kk)\kk} > 0$ (resp. $B_{s(\el)\el} > 0$) for all $\kk$ (resp. $\el$), the following sum $\Sigma$ reduces to
\begin{equation*}
  \Sigma = \sum_{\kk \kk'} \sum_{\el \el'} A_{\kk \kk'} B_{\el \el'} d_{\kk\el\kk'\el'} = \sum_{\kk} \sum_{\el} A_{s(\kk)\kk} B_{t(\el)\el} d_{s(\kk)t(\el) \kk\el}
\end{equation*}
$\Sigma$ is null and sum of positive components, each component is null. However, all $A_{s(\kk)\kk}$ and $B_{t(\el)\el}$ are not null, so that for all $(\kk, \el)$, $d_{s(\kk)t(\el) \kk\el} = 0$ and $\bar{\etaa}_{\kk\el} = \etaa_{s(\kk) t(\el)}$.\\
Now, if $A$ is not a permutation matrix while $\Sigma = 0$ (the same reasoning holds for $B$ or both). Then $A$ owns a column $k$ that contains two non zero elements, say $A_{\kk_1 \kk}$ and $A_{\kk_2 \kk}$. Let $\el \in \{1\dots\m\}$, there exists by assumption $\el'$ such that $B_{\el \el'} \neq 0$. As $\Sigma=0$, both products $A_{\kk_1 \kk} B_{\el \el'} d_{\kk_1 \el \kk \el'}$ and $A_{\kk_2 \kk} B_{\el \el'} d_{\kk_2 \el \kk \el'}$ are zero.
\begin{equation*}
  \begin{cases}
    A_{\kk_1 \kk} B_{\el \el'} d_{\kk_1 \el \kk \el'} & = 0 \\
    A_{\kk_2 \kk} B_{\el \el'} d_{\kk_2 \el \kk \el'} & = 0 \\
  \end{cases}
  \Leftrightarrow
  \begin{cases}
    d_{\kk_1 \el \kk \el'} & = 0 \\
    d_{\kk_2 \el \kk \el'} & = 0 \\
  \end{cases}
  \Leftrightarrow
  \begin{cases}
    \etaa_{\kk_1 \el} = \bar{\etaa}_{\kk \el'} & \\
    \etaa_{\kk_2 \el} = \bar{\etaa}_{\kk \el'} & \\
  \end{cases}
  \Leftrightarrow
  \etaa_{\kk_1 \el} = \etaa_{\kk_2 \el}
\end{equation*}
The previous equality is true for all $\el$, thus rows $\kk_1$ and $\kk_2$ of $\etaa$ are identical, and contradict the assumptions. \proofend 
%
%La réciproque est immédiate.
%\proofend 

%% Lemma for log likelihood ratio of assignments
%\subsection{Likelihood ratio of assignments}

\begin{lemme}\label{lemme:marginalprobabilties}\textbf{}

Let $\mcZ_1$ be the subset of $\mcZ$ of $c$-regular configurations, as defined in Definition~\ref{def:regular}. 
Let $\mathbb{S}^{\g} = \{ \bpi = (\pii_1, \pii_2, \dots, \pii_\g) \in [0, 1]^\g : \sum_{k=1}^\g \pii_k = 1\}$ be the $g$-dimensional simplex and
denote $\mathbb{S}^{\g}_c = \mathbb{S}^{\g} \cap [c, 1-c]^\g$. Then there exists two positive constants $M_c$ and $M'_c$ such that for all $\bz$, $\bzs$ in $\mcZ_1$ and all $\bpi \in \mathbb{S}_c^{g}$
\begin{eqnarray*}
\left| \log \prob(\bz; \widehat{\bpi}(\bz)) - \log \prob(\bzs; \widehat{\bpi}(\bzs)) \right| & \leq & M_c \|\bz - \bzs \|_0 
\end{eqnarray*}
\end{lemme}

\proofbegin
Consider the entropy map $H: \mathbb{S}^{g} \to \mathbb{R}$ 
defined as $H(\bpi) = - \sum_{k = 1}^{\g} \pii_k \log(\pii_k)$. The gradient $\nabla H$ is uniformly bounded by $\frac{M_c}{2} = \log \frac{1-c}{c}$ in $\|.\|_{\infty}$-norm over $\mathbb{S}^{\g} \cap [c, 1-c]^\g$. 
Therefore, for all $\bpi$, $\bpis \in \mathbb{S}^{\g} \cap [c, 1-c]^\g$, we have 
\begin{equation*}
| H(\bpi) - H(\bpis) | \leq \frac{M_c}{2} \| \bpi -\bpis \|_1
\end{equation*}
To prove the  inequality, we remark that $\bz \in \mcZ_1$ translates to $\widehat{\bpi}(\bz) \in \mathbb{S}^{\g} \cap [c, 1-c]^\g$, 
that $\log \prob(\bz; \widehat{\bpi}(\bz)) - \log \prob(\bzs; \widehat{\bpi}(\bzs)) = n [ H(\widehat{\bpi}(\bz)) - H(\widehat{\bpi}(\bzs)) ]$ 
and finally that $\|\widehat{\bpi}(\bz) - \widehat{\bpi}(\bzs) \|_1 \leq \frac{2}{n} \| \bz - \bzs \|_0$. 
\proofend

\begin{corollaire}\label{cor:marginalprobabilties}\textbf{}
Let $\bzs$ (resp. $\bws$) be $c/2$-regular and $\bz$ (resp. $\bw$) at $\|.\|_0$-distance $c/4$ of $\bzs$ (resp. $\bws$). Then, for all $\btheta \in \bTheta$
\begin{equation*}
\log \frac{\prob(\bz, \bw; \btheta)}{\prob(\bzs, \bws; \bthetas)} \leq \bigO_P(1) \exp \left\{ M_{c/4} ( \|\bz - \bzs \|_0 + \|\bw - \bws \|_0 ) \right\} \\ 
% & \times \exp \left\{ M'_{c/2} ( \n \|\bpi(\bzs) - \bpis \|_1 + \dd \|\brho(\bws) - \brhos \|_1 ) \right\}
\end{equation*}
\end{corollaire}

\proofbegin
Note then that:
\begin{align*}
\frac{\prob(\bz, \bw; \btheta)}{\prob(\bzs, \bws; \bthetas)} & = & \frac{\prob(\bz, \bw; \bpi, \brho)}{\prob(\bzs, \bws; \bpis, \brhos)} = 
\frac{\prob(\bz, \bw; \bpi, \brho)}{\prob(\bzs, \bws; \widehat{\bpi}(\bzs), \widehat{\brho}(\bws))} \frac{\prob(\bzs, \bws; \widehat{\bpi}(\bzs), \widehat{\brho}(\bws))}{\prob(\bzs, \bws; \bpis, \brhos)} \\
& \leq & \frac{\prob(\bz, \bw; \widehat{\bpi}(\bz), \widehat{\brho}(\bw))}{\prob(\bzs, \bws; \widehat{\bpi}(\bzs), \widehat{\brho}(\bws))} \frac{\prob(\bzs, \bws; \widehat{\bpi}(\bzs), \widehat{\brho}(\bws))}{\prob(\bzs, \bws; \bpis, \brhos)} \\
& \leq & \exp \left\{ M_{c/4} ( \|\bz - \bzs \|_0 + \|\bw - \bws \|_0 ) \right\} \times \frac{\prob(\bzs, \bws; \widehat{\bpi}(\bzs), \widehat{\brho}(\bws))}{\prob(\bzs, \bws; \bpis, \brhos)} \\
% & \leq & \exp \left\{ M_{c/4} ( \|\bz - \bzs \|_0 + \|\bw - \bws \|_0 ) \right\} \\ 
% & & \times \exp \left\{ M'_{c/2} ( \n \|\bpi(\bzs) - \bpis \|_1 + \dd \|\brho(\bws) - \brhos \|_1 ) \right\}
& \leq & \bigO_P(1) \exp \left\{ M_{c/4} ( \|\bz - \bzs \|_0 + \|\bw - \bws \|_0 ) \right\}  \\
\end{align*}
where the first inequality comes from the definition of $\widehat{\bpi}(\bz)$ and $\widehat{\brho}(\bw)$ and the second from Lemma~\ref{lemme:marginalprobabilties} and the fact that $\bzs$ and $\bz$ (resp. $\bws$ and $\bw$) are $c/4$-regular. 
% \textcolor{red}{and the last from Lemma~\ref{lemme:marginalprobabilties} and the fact that $\widehat{\bpi}(\bzs)$ and $\bpis$ (resp. $\widehat{\brho}(\bws)$ and $\brhos$) are in $\mathbb{S}^{\g}_{c/2}$ (resp. $\mathbb{S}^{\m}_{c/2}$).}
Finally, local asymptotic normality of the MLE for multinomial proportions ensures that $\frac{\prob(\bzs, \bws; \widehat{\bpi}(\bzs), \widehat{\brho}(\bws))}{\prob(\bzs, \bws; \bpis, \brhos)} = \bigO_P(1)$. 
\proofend

%%%%%%%%%%
%\section{Proofs regarding the estimators}